\documentclass[11pt,a4paper]{amsart}
\usepackage{amsmath,amsfonts,amssymb,epsf,mathrsfs}
\usepackage[T1]{fontenc}

\usepackage{graphicx,diagrams}
\usepackage[backref=page,colorlinks,bookmarksopen=true]{hyperref}

\let\bbbibitem\bibitem
\renewcommand{\bibitem}[2][]{\bbbibitem[#1]{#2}\label{#2}}

\def\fin{\hfill\hbox{\hskip .2cm $\square$}\medskip}

\newtheorem{theo}{Theorem}[section]
\newtheorem{lemma}[theo]{Lemma}
\newtheorem{prop}[theo]{Proposition}

\newenvironment{demo}[1][\hspace{-3pt}]{{\noindent\em Proof #1.~ }}{\fin}

\def\cal{\mathcal}
\def\a{\alpha}
\def\b{\beta}
\def\deg{{\rm deg}\,}
\def\g{\gamma}
\def\G{\Gamma}
\def\R{{\mathbb R}} 
\def\d{{\rm d}}
\def\Ad{{\rm Ad}}
\def\ad{{\rm ad}}

\def\PU{{\rm PU}}
\def\SU{{\rm SU}}
\def\GL{{\rm GL}}
\def\End{{\rm End}}
\def\Hom{{\rm Hom}}
\def\U{{\rm U}}
\def\SL{{\rm SL}}
\def\SO{{\rm SO}}

\def\Ker{{\rm Ker\,}}
\def\Im{{\rm Im\,}}
\def\C{{\mathbb C}}
\def\B{{\mathbb H}_{\mathbb C}}
\def\Bm{{\mathbb H}_{\mathbb C}^m}

\def\N{{\mathbb N}}
\def\Z{{\mathbb Z}}
\def\R{{\mathbb R}}
\def\fd{\longrightarrow}

\def\la{\langle}
\def\ra{\rangle}

\def\om{\omega}

\def\rk{{\rm rk}\,}

\def\t{\theta}
\def\ds{\displaystyle}

\def\Om{\Omega^{1}}
\def\o{\omega}
\def\T{T^{1}}
\def\X{{\cal X}}

\def\U{{\rm U}}
\def\gg{{\mathfrak g}}
\def\kk{{\mathfrak k}}

\def\pp{{\mathfrak p}}
\def\r{\rho}
\def\dbar{\bar\partial}

\def\F{{\cal F}}

\def\got{\mathfrak}
\def\tr{{\rm tr}}
\def\trans{^t\hskip -1pt}

\def\V{{\mathbb V}}
\def\W{{\mathbb W}}
\def\E{{\mathbb E}}

\begin{document}

\title[Representations of complex hyperbolic lattices]{Representations of complex hyperbolic lattices \linebreak into rank 2 classical Lie Groups of Hermitian type}

\author{Vincent Koziarz} 
\author{Julien Maubon}
\address{Institut Elie Cartan, Université
    Henri Poincaré, B. P. 239, F-54506 Vand\oe uvre-lès-Nancy Cedex,
    France}
\email{koziarz@iecn.u-nancy.fr}  
\email{maubon@iecn.u-nancy.fr}

\date{March 2007}

\sloppy

\begin{abstract}
Let $\G$ be a torsion-free uniform lattice of $\SU(m,1)$, $m>1$. Let $G$ be either $\SU(p,2)$ with $p\geq 2$, ${\rm Sp}(2,\R)$ or $\SO(p,2)$ with $p\geq 3$. The symmetric spaces associated to these $G$'s are the classical bounded symmetric domains of rank 2, with the exceptions of $\SO^\star(8)/\U(4)$ and $\SO^\star(10)/\U(5)$. Using the correspondence between representations of fundamental groups of Kähler manifolds and Higgs bundles we study representations of the lattice $\G$ into $G$. We prove that the Toledo invariant associated to such a representation satisfies a Milnor-Wood type inequality and that in case of equality necessarily $G=\SU(p,2)$ with $p\geq 2m$ and the representation is reductive, faithful, discrete, and stabilizes a copy of complex hyperbolic space (of maximal possible induced holomorphic sectional curvature) holomorphically and totally geodesically embedded in the Hermitian symmetric space $\SU(p,2)/{\rm S}(\U(p)\times\U(2))$, on which it acts cocompactly.       
\end{abstract}

\maketitle

\section{Introduction}

Let $\G$ be a (torsion free) uniform lattice in the Lie group $\SU(m,1)$. We are interested here in representations, i.e. homomorphisms, of $\G$ in a Lie group $G$ of Hermitian type, that is a connected semisimple Lie group with finite center and no compact factor whose associated symmetric space $\X=G/K$ is Hermitian ($K$ is a maximal compact subgroup of $G$). We will always assume that $G$, hence $\X$, are irreducible. The classical groups of Hermitian type are $\SU(p,q)$, ${\rm Sp}(n,\R)$, $\SO^\star(2n)$ and $\SO_0(n,2)$ whose associated symmetric spaces' real ranks are respectively ${\rm min}(p,q)$, $n$, $[n/2]$, and 2.

The Toledo invariant is a number naturally associated to such a representation $\rho:\G\fd G$, and it has been recognized over the years to play a fundamental role.  It is defined as follows. Let $\Bm=\SU(m,1)/{\rm S}(\U(m)\times\U(1))$ be complex hyperbolic $m$-space, and let $f$ be any (smooth) $\rho$-equivariant map $\Bm\fd\X$. The symmetric space $\X$ being Hermitian, it is a Kähler manifold, and its Kähler form $\o_\X$ may be pulled-back by $f$ to give a 2-form $f^\star\o_\X$ on $\Bm$ which goes down by $\rho$-equivariance to a form on the closed complex hyperbolic manifold $M=\G\backslash\Bm$. We will make no difference between $\G$-invariant objects on $\Bm$ and the corresponding objects on $M$. For example, $f^\star\o_\X$ is either a 2-form on $\Bm$ or a 2-form on $M$, depending on the context. Similarly, we denote by $g$ and $\o$ the invariant metric and Kähler form of $\Bm$, as well as the induced metric and Kähler form on $M$.  
Now, the de Rham cohomology class in ${\rm H}^2_{dR}(M)$ defined by $f^\star\o_\X$ depends only on $\rho$, not on $f$, and will be denoted by $[\rho^\star\o_\X]$. This class is then evaluated against the Kähler class of $M$ to give the desired number
$$
\tau(\rho)=\frac{1}{2m}\int_M\la\rho^\star\o_\X,\o\ra\, dV=\frac{1}{m!}\int_M\rho^\star\om_\X\wedge\o^{m-1}
$$
where $\rho^\star\om_\X$ is any representative of $[\rho^\star\om_\X]$, $\la.,.\ra$ is the scalar product induced by $g$ on 2-forms and $dV=\frac{1}{m!}\,\o^m$ is the Riemannian volume form of $\Bm$ (or $M$).
          
The Toledo invariant is of particular interest because

(1) it is constant on connected components of the space of representations $\Hom(\G,G)$;

(2) it satisfies a Milnor-Wood type inequality, namely $\tau(\rho)$ is bounded in absolute value by a quantity depending only on the (real) rank of the symmetric space $\X$ and the volume of $M=\G\backslash\Bm$. More precisely, if the Riemannian metrics on $\Bm$ and $\X$ are normalized so that the minimum of their holomorphic sectional curvatures is $-1$ (so that the holomorphic sectional curvature of $\X$ is pinched between $-1$ and $-1/\rk\X$), the following holds:
$$
|\tau(\rho)| \leq \rk\X\, {\rm Vol}(M).
$$

(3) maximal representations, i.e. representations $\rho$ for which $|\tau(\rho)|= \rk\X\, {\rm Vol}(M)$, are expected to be of a very special kind, and therefore rigidity results should follow.

The Toledo invariant appeared for the first time in 1979 Toledo's paper~\cite{T1} and more explicitly in~\cite{T2}, where the Milnor-Wood inequality (2) was proved for $m=1$ and $\rk\X=1$, namely $\G$ is the fundamental group of a Riemann surface and $G=\SU(n,1)$. Toledo also proved that maximal representations are faithful with discrete image, and stabilize a complex line in the complex hyperbolic $n$-space $\X$, generalizing Goldman's results for $G={\rm SL}(2,\R)$~\cite{G1,G2}. At approximately the same time, Corlette, using a very similar invariant, the volume of the representation, established in~\cite{C} the same kind of result for $m\geq 2$ and $G=\SU(n,1)$. An immediate corollary is that a uniform lattice in $\SU(m,1)$ can not be deformed non-trivially in $\SU(n,1)$, $n\geq m\geq 1$, a result first obtained by Goldman and Millson in~\cite{GM} using different methods. These results have been extended to the non-uniform case in~\cite{BI1,KM} (the definition of the Toledo invariant must be modified). Therefore the case where the rank of the symmetric space associated to $G$ is 1 is now settled.  

Using the work of Domic-Toledo~\cite{DT} and  Clerc-{\O}rsted~\cite{CO}, Burger and Iozzi obtained in~\cite{BI1} the Milnor-Wood inequality (2) in full generality. Since then, much progress have been made  and maximal representations of the fundamental group of a Riemann surface into such general groups of Hermitian type are well understood~\cite{BI2,BI3,BGPG1,BGPG2}.

So far, the problem of representing higher dimensional complex hyperbolic lattices in Lie groups of Hermitian type of rank at least two seems to have been left aside. This is the question we (begin to) address in this paper. Since the Milnor-Wood inequality (2) is known, one would a priori like to focus on point (3) and give a complete description of the maximal representations. Our strategy will be different, and to explain it we need 
to say a word about the available methods to prove the results we mentioned. Essentially, there are two different ways of attacking the problem. The one used by Burger and Iozzi relies on bounded cohomology theory and allows to prove the bound (2) in great generality but gives relatively few informations on the maximal representations, so that a separate study has to be made. The second one, used in~\cite{C,KM} and in~\cite{BGPG1,BGPG2}, relies on harmonic maps and/or Higgs bundles machinery and belongs more to the world of complex differential geometry. Following this approach, as we shall, makes it quite difficult to prove the Milnor-Wood inequality (and in fact no such proof is known in the general case) but once it is proved (in some special cases), maximal representations are easier to understand.       

The Higgs bundle theory was developed by Hitchin~\cite{H1,H2} for Riemann surfaces and Simpson~\cite{S1,S2,S3,S4} in higher dimensions to study fundamental groups of Kähler manifolds and their linear representations. It consists in establishing a correspondence between representations of the fundamental group $\G$ of a Kähler manifold $M$ and purely holomorphic objects, called polystable Higgs bundles over $M$, and then working on these objects with tools from complex geometry. This is in fact only possible for reductive representations, since the construction of Higgs bundles requires harmonic maps, but this restriction will not be a serious issue for our purposes.  
We will explain this correspondence in some details in section~\ref{higgs}. For now let us simply say that to a reductive representation of $\G$ is associated a Higgs bundle on $M$, that is a holomorphic vector bundle $(E,\dbar)$ on $M$ together with a holomorphic (1,0)-form $\t$, the Higgs field, which takes values in (a subbundle of) the bundle of endomorphisms of $E$ and satisfies $[\t,\t]=0$. We shall often write $\t:E\fd E\otimes\Om$, where $\Om$ is the sheaf (of germs) of holomorphic 1-forms on $M$, which we identify with the holomorphic cotangent bundle of $M$.   Polystability is a condition that relates the slope of proper $\t$-invariant saturated subsheaves of $E$ to the slope of $E$ itself (which will be 0 here). Recall that the slope of a saturated sheaf is defined as the quotient of its degree by its rank.    

The Toledo invariant may be defined as before for representations $\rho$ of the fundamental group $\G$ of any closed Kähler manifold $M$ in a Lie group of Hermitian type $G$ ($f$ is then simply a $\rho$-equivariant map from the universal cover ${\widetilde M}$ of $M$ to $\X=G/K$), but in the Higgs bundles setting it is best interpreted as the degree of a complex vector bundle over $M$ as follows. The Hermitian symmetric space $\X$ is Kähler-Einstein and therefore $\o_\X$ is up to a constant the first Chern class of the holomorphic tangent bundle $T\X$ of $\X$. Therefore $f^\star\o_\X$ is up to a constant the first Chern class of the induced bundle $f^\star T\X$ over $M$ so that $\tau(\rho)$ is, again up to a constant, simply the degree of this bundle. To be more precise, if the Ricci curvature tensor of $\X$ is $\lambda_\X\,g_\X$ (remember that the Riemannian metric of $\X$ is normalized to have minimal holomorphic sectional curvature $-1$), we have $\tau(\rho)=-\frac{2\pi}{m!\lambda_\X}\,\deg f^\star T\X$. We remark that if $M$ is complex hyperbolic of complex dimension $m$, ${\rm Vol}(M) = \frac{4\pi}{m! (m+1)}\,\deg \Om$. Hence the Milnor-Wood type bound (2) can be written:
\begin{equation}
\left|\frac{\deg f^\star T\X}{2\lambda_\X}\right|\leq \rk\X\,\frac{\deg\Om}{m+1}\,,\tag{\mbox{2'}} 
\end{equation}
an interpretation we will use constantly. 
If the $\rho$-equivariant map $f:{\widetilde M}\fd\X$ is chosen to be harmonic, the bundle $f^\star T\X\fd M$ is constructed from holomorphic subbundles of the Higgs bundle $E\fd M$ associated to the representation $\rho$. Moreover, the fact that $\rho$ is not valued in the full general linear group but in a smaller group $G$ of Hermitian type implies that the Higgs bundle $E$ has a special structure. The idea is then that Inequality (2') will follow from this particular structure and the polystability condition. This implementation of the Higgs bundle theory has been carried out for Riemann surfaces, for example by Xia~\cite{X} in some special cases, and more generally by Bradlow, Garcia-Prada and Gothen in~\cite{BGPG1,BGPG2}. The Kähler manifold $M$ being a complex curve in their situation makes it quite easy to deduce (2') from the structure of the Higgs bundles and maximal representations can be studied in great details. It is for example possible to count the number of connected components of the moduli space of maximal representations. The reader should consult the papers~\cite{BGPG1,BGPG2} to see the strength of the method in this case.  

When $M=\G\backslash\Bm$ is a (closed) complex hyperbolic manifold of dimension $m\geq 2$, one expects maximal representations to be extremely restricted. In fact they should all be induced by special holomorphic or antiholomorphic totally geodesic embeddings of complex hyperbolic $m$-space into the Hermitian symmetric space $\X$ associated to $G$. 

Inequality (2) or (2') for such higher dimensional lattices turns out to be surprisingly difficult to prove and disappointingly we have been obliged to restrict ourselves to the case where the rank of the symmetric space $\X$ is 2. Our main result is    

\begin{theo}
Let $\G$ be a (torsion-free) uniform lattice in $\SU(m,1)$. Let $G$ be either $\SU(p,2)$ with $p>1$, $\SO_0(p,2)$ with $p\geq 3$ or ${\rm Sp}(2,\R)$.  Finally let $\rho:\G\fd G$ be a representation. 

Then $|\tau(\rho)|\leq 2 {\rm Vol}(\G\backslash\Bm)$. If $m>1$ and $\rho$ is maximal, namely if $|\tau(\rho)|= 2 {\rm Vol}(\G\backslash\Bm)$, then $G=\SU(p,2)$ with $p\geq 2m$, $\rho$ is reductive, faithful, discrete, and stabilizes a holomorphic totally geodesic copy of complex hyperbolic $m$-space of holomorphic sectional curvature $-1/2$ in the Hermitian symmetric space $\X=\SU(p,2)/{\rm S}(\U(p)\times\U(2))$. Moreover $\G$ acts cocompactly on this copy of complex hyperbolic space.   
\end{theo}

Although Higgs bundles are associated to reductive representations, we do not assume that $\rho$ is reductive in the theorem. This is because every representation can be deformed to a reductive one, an operation that does not change the value of the Toledo invariant. Moreover we shall see in section~\ref{nonreductive} that non reductive representations can not be maximal.   

The theorem covers all but two of the classical Lie groups of Hermitian type $G$ whose associated symmetric spaces' rank is 2. The missing ones are $\SO^\star(8)$ and $\SO^\star(10)$. In fact, $\SO^\star(8)$ and $\SO_0(6,2)$ are isogenous and have the same associated symmetric space. This means that if a representation $\G\fd \SO^\star(8)$ lifts to 
${\rm Spin}(6,2)$, then projecting down to $\SO_0(6,2)$ gives the result for this representation as well.
  
As we said, the Milnor-Wood type inequality  $|\tau(\rho)|\leq 2 {\rm Vol}(M)$ is not new, only the proof is. It is given in Section~\ref{sup2} for $G=\SU(p,2)$ and in Section~\ref{sop2} for $G=\SO_0(p,2)$. The case of ${\rm Sp}(2,\R)$ follows from the case of $\SU(2,2)$ since  ${\rm Sp}(2,\R)\subset\SU(2,2)$. 

The theorem in particular says that for $m>1$ there is no maximal representations of a uniform lattice $\G$ of $\SU(m,1)$ in $\SO_0(p,2)$, ${\rm Sp}(2,\R)$ or $\SU(p,2)$ with $p<2m$. Our method indeed yields explicit better bounds on the Toledo invariant in these cases. For representations $\rho: \G\fd\SU(p,2)$, the arguments of~\cite{VZ05} can be adapted to give the following, which is stronger than the Milnor-Wood inequality (2) exactly when $p<2m$:

\begin{prop}\label{otherbound}
Let $\G$ be a (torsion-free) uniform lattice in $\SU(m,1)$ and let $\rho:\G\fd\SU(p,2)$, $p\geq 2$, be a representation. Then $|\tau(\rho)|\leq\frac{2p}{p+2}\frac{m+1}{m}\,{\rm Vol}(\G\backslash\Bm)$.
\end{prop}

If $\rho:\G\fd G$ is a maximal representation, for which as we said $G=\SU(p,2)$ with $p\geq 2m$, we will prove that there exists a maximal holomorphic or antiholomorphic totally geodesic $\rho$-equivariant embedding $\Bm\fd\X=\SU(p,2)/{\rm S}(\U(p)\times\U(2))$, from which the assertions of our main theorem follow. By a maximal embedding $\Bm\fd\X$ we mean an embedding whose image's induced holomorphic sectional curvature is everywhere the greatest possible, namely $-1/2$ with our normalization. See section~\ref{maxemb} for a discussion of these embeddings and a description of the stabilizer in $\SU(p,2)$ of their images in $\X$.  If $f:\Bm\fd\X$ is such a $\rho$-equivariant maximal embedding, we will loosely say that $\rho$ is induced by $f$, although $f$ determines $\rho(\g)$ for $\g\in\G$ only up to composition with an element of $\SU(p,2)$ fixing pointwisely the image of $f$ in $\X$.        

The paper is organized as follows. Section~\ref{higgs} is an overview on how Higgs bundles are constructed from representations of the fundamental group of a Kähler manifold. We say a few words about the corresponding moduli space, the $\C^\star$-action it comes with, and the systems of Hodge bundles that are obtained as fixed points of this action. Section~\ref{sup2} is devoted to the proof of the main theorem when the representation takes values in $\SU(p,2)$, which is the most interesting case. The first subsection is expository, we give there the necessary background on the geometry of the associated Hermitian symmetric space. This is used in the next subsection to describe the particular structure of the Higgs bundles associated to such a representation. The third subsection contains the proof of the Milnor-Wood type inequality and the forth deals with maximal representations. In the fifth we prove Proposition~\ref{otherbound}. Finally, Section~\ref{sop2} follows the lines of Section~\ref{sup2} in the case of $\SO_0(p,2)$: the first subsection describes the associated symmetric space whereas the second is devoted to the Higgs bundles arising in this case and to the proof of the Milnor-Wood type inequality.      

\vspace{0.5cm}

{\em Acknowledgments.} We would like to thank Jean-Louis Clerc and Andrei Teleman for useful discussions and their interest in our work.

\section{Representations of the fundamental group, flat bundles and Higgs bundles}\label{higgs}

In this section we give a short presentation of the links between representations of the fundamental group $\G$ of a Kähler manifold $M$ and Higgs bundles on $M$. To be a little more precise, we will explain in some details how $G$-Higgs bundles are constructed from reductive representations of $\G$ into a linear group $G$. There is in fact a much deeper correspondence (a generalized Hitchin-Kobayashi correspondence) between the moduli space of reductive representations and the moduli space of $G$-Higgs bundles over $M$ with some stability properties. However, we shall not need the full strength of this correspondence (the easy direction suffices), and we will stick to what matters for our purposes. We refer to the original papers of Simpson~\cite{S1,S2,S3,S4} and to~\cite{BGPG1,BGPG2} for details. Our exposition owes a lot to~\cite{BGPG3}.

Let $M$ be a compact Kähler manifold, $\G$ its fundamental group, and ${\widetilde M}$ its universal cover, so that $M=\G\backslash{\widetilde M}$. Let $G$ be a real connected semisimple Lie group with finite center and no compact factor and $K$ a maximal compact subgroup of $G$. Finally, let $\rho$ be a representation $\G\fd G$.    

\subsection{Real Higgs equations}\hfill

Let $P_G$ be the flat principal $G$-bundle ${\widetilde M}\times_\r G$ on $M$ associated to the representation $\rho$.  

A metric on $P_G$ is a reduction of the structure group $G$ of $P_G$ to its maximal compact subgroup $K$, namely, a $K$-principal subbundle $P_K$ of $P_G$. This is the same thing as a section of the associated bundle $P_G\times_G \X\simeq P_G/K$ over $M$. In our setting, since $P_G$ is flat, this associated bundle is isomorphic to ${\widetilde M}\times_\r \X$, and a section of this bundle is given by a $\r$-equivariant map $f:{\widetilde M}\fd\X$. In this case, the $K$-principal bundle $G\fd \X=G/K$ can be pulled-back by $f$ to give a $K$-principal bundle $f^\star G\subset {\widetilde M}\times G$ over ${\widetilde M}$. This bundle goes down under the action of $\G$ and yields the $K$-principal bundle $P_K\subset P_G$ over $M$. Note that $P_G$ is recovered as the bundle $P_K\times_K G$ associated to $P_K$ via the action of $K$ on $G$ by left translations.   

\begin{diagram} 
              &            & f^\star G \subset {\widetilde M}\times G &        & G          \\
              & \ldTo      & \dTo                         &        & \dTo    \\
P_K \subset P_G &            & {\widetilde M}                           & \rTo^f & \X \\
\dTo          & \ldTo      &            &        &        \\
 M            &            &            &        &        \\
\end{diagram}

Let $\tilde\o_G$ be the flat connection 1-form on the trivial $G$-bundle ${\widetilde M}\times G\fd {\widetilde M}$ : if $X\in T{\widetilde M}$ and $A^*$ is the left invariant vector field on $G$ corresponding to $A\in\gg$, $\tilde\o_G(X,A^*)=A$.  This form goes down under the $\G$-action to give the flat connection $\o_G$ on $P_G$. On the bundle $G\fd\X$ we have the usual invariant connection $\lambda$ defined by $\lambda(A^*)=A_\kk$, where $A_\kk$ is the $\kk$-component of $A\in\gg$ in the Cartan decomposition $\gg=\kk\oplus\pp$. Let $\tilde\o_K=f^\star\lambda$ be the induced connection 1-form on the pull-back $f^\star G\fd {\widetilde M}$. Again, $\tilde\o_K$ is $\G$-invariant and gives a connection 1-form $\o_K$ on $P_K\fd M$. 

For $X\in T_x{\widetilde M}$, we have $\d_xf(X)\in T_{f(x)}\X=T_{gK}G/K=g_\star T_{eK}G/K=g_\star\pp$.
Hence we can define a form $\tilde \Theta\in \Om({\widetilde M}\times G)\otimes\pp$ by $\tilde\Theta(X,A^*)=g_\star^{-1}\d f(X)$. If we restrict it to $f^\star G$, we have that $\tilde\Theta(X,A^*)=A_\pp$. Hence, on $f^\star G$, we have $\tilde\o_G=\tilde\o_K+\tilde\Theta$. $\Theta$ is $\G$-invariant and gives an element $\Theta$ of $\Om(P_K)\otimes\pp$ so that 
$\o_G=\o_K+\Theta$ on $P_K$. The form $\Theta$ behaves well under the right action of $K$ on $P_K$: $R_k^\star\Theta=\Ad(k^{-1})\Theta$. Moreover $\Theta$ vanishes on vectors tangent to the fibers of $P_K\fd M$. Hence $\Theta$ can be seen as a 1-form on $M$ with values in the vector bundle $P_K\times_\Ad\pp\fd M$ associated to $P_K$ via the adjoint action of $K$ on $\pp$. One should remark that this vector bundle is nothing but the quotient under $\G$ of the pull-back $f^\star T\X$ of the tangent bundle $T\X$ of $\X$.            

Let $\d_G$, $F_G$ and $\d_K$, $F_K$ be the covariant exterior derivatives and the curvature forms of the connection 1-forms $\o_G$ and $\o_K$. We have for example $\d_G=\d+\ad(\o_G)$ and $F_G=\d\o_G+\frac{1}{2}[\o_G,\o_G]$. 

Since $\o_G=\o_K+\Theta$ is flat, we have $0=F_G=F_K+\frac{1}{2}[\Theta,\Theta]+\d_K\Theta$. Decomposing according to $\gg=\kk\oplus\pp$, we obtain 
$$
\left\{
\begin{array}{r}
F_K+\frac{1}{2}[\Theta,\Theta]=0\\
\d_K\Theta=0
\end{array}
\right.
$$   

Moreover, Corlette proved in \cite{C} that if the representation $\r$ is reductive, that is if the Zariski closure of $\r(\G)$ is a reductive subgroup of $G$, then there exists a harmonic $\r$-equivariant map $f:{\widetilde M}\fd\X$. In our setting, this means that $\Theta$ can be chosen to satisfy the additional condition 
$$
\d_K^\star\Theta=0.
$$  

Summing up, we see that a reductive representation $\r:\G\fd G$ is equivalent to the data of a connection $\d_K$ on a $K$-principal bundle $P_K\fd M$ and an element $\Theta\in\Om(M,P_K\times_\Ad\pp)$ satisfying 
$$
\left\{
\begin{array}{r}
F_K+\frac{1}{2}[\Theta,\Theta]=0\\
\d_K\Theta=0\\
\d_K^\star\Theta=0
\end{array}
\right.\mbox{ (real Higgs equations) }
$$

So far we have not taken into account the fact that $M$ is a Kähler manifold. 

\subsection{Complex Higgs equations}\hfill

Let $G_\C$, $K_\C$ and $\gg_\C=\kk_\C\oplus \pp_\C$ be the complexifications of $G$, $K$, and $\gg$. 

The $K$-principal bundle $P_K$ can be extended to a $K_\C$-principal bundle $P_{K_\C}=P_K\times_K K_\C$. We still denote by $\d_K$ the extension of the covariant exterior derivative to $P_{K_\C}$. $\Theta$ can be extended by $\C$-linearity to an element of $\Om_\C(M,P_{K_\C}\times_\Ad\pp_\C)=\Omega^0(T^\star_\C M\otimes (P_{K_\C}\times_\Ad\pp_\C))$. Decomposing according to types we write $\d_K=\d_K^{1,0}+\d_K^{0,1}$ and $\Theta=\Theta^{1,0}+\Theta^{0,1}$. We set $\t=\Theta^{1,0}\in\Omega^{1,0}(M,P_{K_\C}\times_\Ad\pp_\C)$. If we call $\tau$ the involution of $\gg$ defined by $\tau(A)=-^{t}\bar A$, then $\Theta=\t-\tau(\t)$, for $\t$ comes from a real form.

It has been shown by Sampson~\cite{Sa} that the harmonicity of the map $f$ and the fact that $\X$ has non positive complexified sectional curvature implies that $f$ is pluriharmonic, which can be written $\d_K^{0,1}\t=0$. Moreover the complexified sectional curvature of $\X$ vanishes on the image of the (1,0)-part of $f$, i.e. on the image of $\t$. Since here $\X$ is a symmetric space, the curvature is given by the Lie bracket and the vanishing of the complexified sectional curvature just means that $[\t,\t]=0$. From this it follows that ${(\d_K^{0,1})}^2=0$, namely that $\d_K^{0,1}$ defines a holomorphic structure on the principal bundle $P_{K_\C}\fd M$ and on the associated complex vector bundle $P_{K_\C}\times_\Ad\pp_\C$. The pluriharmonicity of $f$ therefore means that $\t$ is holomorphic for this holomorphic structure.  

Looking at the real Higgs equations, we also obtain that $F_K^{1,1}-[\t,\tau(\t)]=0$.

Therefore, starting from a reductive representation $\r$, we have constructed a $K_\C$-principal bundle $P_{K_\C}=P_K\times_K K_\C$ over $M$ equipped with a complex structure $\bar\partial=\d_K^{0,1}$, and a section $\t$ of $\Om\otimes (P_{K_\C}\times_\Ad\pp_\C)$, where $\Om$ denotes the holomorphic 1-forms on $M$, satisfying 

$$
\left\{
\begin{array}{r}
\left[\theta,\theta\right]=0\\
F_K^{1,1}-[\t,\tau(\t)]=0\\
\end{array}
\right.\mbox{ (complex Higgs equations) }
$$      

\subsection{Higgs bundles and stability conditions}\hfill

We will call the data of a $K_\C$-principal bundle $P_{K_\C}$ together with a holomorphic structure $\dbar$ and a section $\t$ of $\Om\otimes (P_{K_\C}\times_\Ad\pp_\C)$ satisfying  $\left[\theta,\theta\right]=0$ a $G$-Higgs principal bundle on $M$. 
This is a purely holomorphic object. The section $\t$ is called a Higgs field. 
   
The remaining differential geometric data, namely the reduction $P_K$ of the structure group of $P_{K_\C}$ to $K$ and the connection $\d_K$ on $P_K$ such that $\bar\partial=\d_K^{0,1}$ and $F_K^{1,1}-[\t,\tau(\t)]=0$, can be rephrased in  holomorphic terms with the help of adapted notions of stability for associated vector bundles. If ${\mathbb E}$ is a complex vector space on which $G_\C$ acts, we can form the associated vector bundle $E=P_{K_\C}\times_{K_\C}{\mathbb E}$. The holomorphic structure on $P_{K_\C}$ induces a holomorphic structure $\dbar$ on $E$. Moreover, $\pp_\C$ can be seen as a subspace of $\End({\mathbb E})$, hence we may consider $\t$ as a holomorphic (1,0)-form on $M$ with values in a subbundle $P$ of the bundle $\End(E)$ (we will write $\t:E\fd E\otimes\Om$). We will call the holomorphic bundle $(E,\dbar)$ together with $\t$ a $G$-Higgs vector bundle. We will often abbreviate $(E,\dbar,\t)$ by $(E,\t)$.

Now the connection $\d_K$ induces a connection on $E$ compatible with its complex structure and the curvature of this connection can be used to compute the degree (and the slope) of saturated subsheaves of $E$. The condition $F_K^{1,1}-[\t,\tau(\t)]=0$ implies that $E$ is polystable in the following sense~\cite{S1}. 

A $G$-Higgs vector bundle $(E,\dbar,\t)$ is called stable (resp. semistable) if for every saturated subsheaf $\F$ of $E$ such that $\t(\F)\subset\F\otimes\Om$ and $0<\rk\F<\rk E$, the slope $\mu(\F)=\frac{\deg\F}{\rk \F}$ of $\F$ is smaller (resp. not bigger) than the slope $\mu(E)=\frac{\deg E}{\rk E}$ of $E$. It is called polystable if it is the sum of stable $G$-Higgs vector bundles of the same slope.  

A subsheaf $\F$ of $E$ such that $\t(\F)\subset\F\otimes\Om$ is called a $\t$-invariant subsheaf or a Higgs subsheaf of $E$. Note that if the Higgs bundle $E$ is polystable and $\F$ is a proper Higgs subsheaf of $E$ with $\mu(\F)=\mu(E)$ then $\F$ is in fact a Higgs subbundle of $E$ and $E$ splits as the direct sum of $\F$ with another Higgs subbundle of the same slope~\cite{S1}.

\subsection{\texorpdfstring{Moduli space, $\C^\star$-action and systems of Hodge bundles}{Moduli space, C*-action and systems of Hodge bundles}}\hfill

Let $P_{K_\C}$ be a fixed $K_\C$-principal bundle on $M$, $\E$ a fixed vector space over $\C$ on which $G_\C$, hence $K_\C$, acts, and let $E=P_{K_\C}\times_{K_\C}\E$ be the associated vector bundle. Consider the space of all holomorphic structures $\dbar$ on $P_{K_\C}\times_{\Ad K_\C}\kk_\C$ and all (1,0)-forms $\t$ on $M$ taking values in the bundle $P=P_{K_\C}\times_{\Ad K_\C}\pp_\C\subset\End E$ such that: 

- $\t$ is holomorphic w.r.t. the complex structure induced on $P$ by $\dbar$;

- the complex structure on $E$ induced by $\dbar$, still denoted by $\dbar$,  turns $(E,\dbar,\t)$ into a polystable $G$-Higgs bundle.   

The group ${\mathscr H}_\C=P_{K_\C}\times_\Ad K_\C$ of gauge transformation of $P_{K_\C}$ acts on this space of polystable $G$-Higgs vector bundles $(E,\dbar,\t)$ by pull-back of the holomorphic structure and conjugacy of the Higgs field. There is a corresponding moduli space ${\cal M}$ (at least if $M$ is a projective variety, see~\cite{S2,S4}), which has the structure of an analytic space. 

A very important feature of this moduli space is that it comes with a natural $\C^\star$-action given by $t.[E,\dbar,\t]=[E,\dbar,t\t]$ for $t\in\C^\star$. Moreover, Simpson proved in~\cite{S4} that for any  $[E,\dbar,\t]\in{\cal M}$, the limit of $[E,\dbar,t\t]$ as $t\in\C^\star$ goes to zero exists and is unique. The limit is therefore a fixed point of the $\C^\star$-action, and this implies that it has the structure of a so-called system of Hodge bundles~\cite{S1,S2}. More precisely, this means that $E$ with the limiting holomorphic structure  splits holomorphically  as a sum $E_1\oplus\ldots\oplus E_k$ of holomorphic vector bundles and that the limiting Higgs field in $P$ is given by a collection of holomorphic maps $\t_i:E_i\fd E_{i+1}\otimes\Om$, $1\leq i \leq k$ (with the convention that $E_{k+1}=\{0\}$). We will abuse notation and use the following kind of diagram:
$$
E_1\stackrel{\t_1}{\fd}E_2\stackrel{\t_2}{\fd}\ldots\stackrel{\t_{k-1}}{\fd}E_k\stackrel{\t_k}{\fd} 0
$$      
to denote such a system of Hodge bundles.

\section{The case $G=\SU(p,2)$}\label{sup2}

In this section we prove the main theorem for representations into $G=\SU(p,2)$. However, some arguments and results are valid in the general case $G=\SU(p,q)$, $p\geq q$, and therefore we will specialize to the case $q=2$ only when necessary.

\subsection{\texorpdfstring{The Hermitian symmetric space $\SU(p,q)/{\rm S}(\U(p)\times\U(q))$}{The Hermitian symmetric space SU(p,q)/S(U(p)xU(q))}}\label{supq}\hfill

\subsubsection{General facts}\hfill

The reader should consult~\cite{Helgason,Satake} for details about this section.

Let $\E$ be a complex vector space of dimension $p+q$, with $p\geq q\geq 1$, endowed with a non-degenerate Hermitian form $F$ of signature $(p,q)$. Let $\W$ be a $q$-dimensional complex subspace of $\E$  on which $F$ is negative-definite, and let $\V$ be its $F$-orthogonal complement, so that $\E=\V\oplus\W$.
The symmetric space $\X$ is defined as the space of all $q$-dimensional complex subspaces of $\mathbb E$ on which $F$ is negative-definite. It is an open submanifold of the complex Grassmannian of $q$-planes of $\mathbb E$. When $q=1$, $\X$ is complex hyperbolic space of (complex) dimension $p$ which we denote by $\B^p$.

After an appropriate choice of basis, we see that the group $G={\SU}(p,q)$ acts transitively on $\X$ by analytic isomorphisms, while the isotropy subgroup $K$ of $G$ at $\W$ is identified with the maximal compact subgroup ${\rm S}(\U(p)\times\U(q))$ of $\SU(p,q)$, so that $\X$ can be identified with ${\SU}(p,q)/{\rm S}(\U(p)\times\U(q))$.

Let $\gg$ be the Lie algebra of $G$, $\kk\subset \gg$ the Lie algebra of $K$
and $\gg=\kk\oplus \pp$ the corresponding Cartan decomposition. We have the following matrix expressions:
$$
\got k=\left\{\left(\begin{array}{cc} X_1 & 0\\ 0 & X_2
\end{array}\right)\ ,\ X_1\in M_p(\C),\ X_2\in M_q(\C),\ \trans\bar X_i=-X_i\ (i=1,2),\ \tr X_1+\tr X_2=0\right\},
$$ 
$$
\got p=\left\{\left(\begin{array}{cc} 0 & A\\ \trans\bar A & 0
\end{array}\right)\ ,\ A\in M_{p,q}(\C)\right\}\,\simeq_\R\,\Hom_\C(\W,\V).
$$ 

The tangent space $T_o\X$ at $o\in\X$ will be identified with $\got p$. More generally, the tangent bundle $T\X$ of $\X$ is the bundle $G\times_{\Ad K}\pp$ associated to the $K$-principal bundle $G\fd \X=G/K$ via the adjoint action of $K$ on $\pp$.  
The complex structure $J$ on $T_o\X$ is given by
$$
J\left(\begin{array}{cc} 0 & A\\ \trans\bar A & 0
\end{array}\right)=\left(\begin{array}{cc} 0 & iA\\ -i\trans\bar A & 0
\end{array}\right)
$$
whereas the $G$-invariant Kähler metric $g_\X$ on $\X$ is defined at $o$ by
$$
g_\X(X,Y)=2\tr\left(YX\right)=4{\rm Re}\,\tr\left(\trans\bar BA\right),\ {\rm if}\ X=\left(\begin{array}{cc} 0 & A\\ \trans\bar A & 0
\end{array}\right),\,Y=\left(\begin{array}{cc} 0 & B\\ \trans\bar B & 0
\end{array}\right)\in \got p.
$$
The corresponding Kähler form will be denoted by $\o_\X=g_\X(J.,.)$.  

The complexifications $G_\C$ of $G$ and $K_\C$ of $K$ are respectively  $\SL(p+q,\C)$ and ${\rm S}(\GL(p,\C)\times\GL(q,\C))$. 
The Lie algebra $\gg_\C$ of $G_\C$ splits as $\kk_\C\oplus \pp_\C$ where $\kk_\C$ is the Lie algebra of $K_\C$ and 
$$
\pp_\C=\pp\otimes\C=\left\{\left(\begin{array}{cc} 0 & A\\ B & 0
\end{array}\right)\ ,\ A\in M_{p,q}(\C)\ ,\ B\in M_{q,p}(\C)\right\}\,\simeq_\C\, \Hom_\C(\W,\V)\oplus\Hom_\C(\V,\W).
$$ 
The extended complex structure $J\otimes{\rm Id}$ acting on $\pp_\C$ has two eigenspaces 
$$
\pp^{1,0}=\left\{\left(\begin{array}{cc} 0 & A\\ 0 & 0\end{array}\right)\ ,\ A\in M_{p,q}(\C)\right\}\,\simeq_\C \Hom_\C(\W,\V)
$$
and 
$$
\pp^{0,1}=\left\{\left(\begin{array}{cc} 0 & 0\\ A & 0\end{array}\right)\ ,\ A\in M_{q,p}(\C)\right\}\,\simeq_\C \Hom_\C(\V,\W).
$$
The complexified tangent bundle $T^\C\X$ of $\X$ is isomorphic to $G\times_{\Ad K}\pp_\C\simeq(G\times_K K_\C)\times_{\Ad K_\C}\pp_\C$, whereas the holomorphic tangent bundle $T^{1,0}\X$ is isomorphic to $G\times_{\Ad K}\pp^{1,0}\simeq(G\times_K K_\C)\times_{\Ad K_\C}\pp^{1,0}$. There is a natural Hermitian metric on the holomorphic tangent bundle of $\X$ given on $T_o^{1,0}\X\simeq M_{p,q}(\C)$ by $h(A,B)=4\tr\left(\trans\bar BA\right)$.   

The holomorphic sectional curvature for the complex line $\la X\ra$ generated by a nonzero $X=\left(\begin{array}{cc} 0 & A\\ \trans\bar A & 0\end{array}\right)\in T_o\X$, or equivalently by a non-zero $\left(\begin{array}{cc} 0 & A\\ 0 & 0
\end{array}\right)\in T_o^{1,0}\X$,
  is given by
$$
K(\la X\ra)=-\frac{\tr\bigl(\left(\trans\bar AA\right)^2\bigr)}{\left(\tr\left(\trans\bar AA\right)\right)^2}.
$$
This formula shows that $K(\la X\ra)$ is pinched between $-1$ and $-1/q$ and that $K(\la X\ra)=-1/q$ if and only if the column vectors of $A$ are pairwise orthogonal and have the same norm (for the standard Hermitian scalar product in $\C^p$).   




The metric $g_\X$ is Einstein and with our normalization, its Ricci curvature tensor is $-\frac{p+q}{2}\,g_\X$. 

\subsubsection{Maximal embeddings}\label{maxemb}\hfill

There is a natural identification of $\X=\SU(p,q)/{\rm S}(\U(p)\times\U(q))$ with the space $\{Z\in M_{p,q}(\C)\ ,\ I_q-{\trans\bar Z}Z>0\}$~\cite{Satake}. Therefore if $m\leq p/q$, we have a holomorphic totally geodesic embedding of complex hyperbolic space $\Bm=\SU(m,1)/{\rm S}(\U(m)\times\U(1))$ into $\X$ given by

\begin{equation}
\Bm\ni z=\left(\begin{array}{c}z_1\\ z_2\\ \vdots \\ z_m
        \end{array}\right)\longmapsto
Z=\left(\begin{array}{c}z_1 I_q\\ z_2 I_q\\ \vdots \\ z_m I_q\\ 0_{p-qm,q}
        \end{array}\right)\in\X.\,\tag{$\star$} 
\end{equation}

This shows that $\X$ contains totally geodesic copies of complex hyperbolic $m$-space of holomorphic sectional curvature $-1/q$, for all $m\leq p/q$. The next lemma implies that $[p/q]$ is the maximal possible dimension of such submanifolds: 

\begin{lemma}\label{tangentmax}
Let $S$ be a complex linear subspace of $T_o^{1,0}\X$. If, for every nonzero $A\in S$, $K(\la A\ra)=-1/q$, then ${\rm dim}_\C S\leq p/q$.
\end{lemma}

\begin{demo}
The metric and the holomorphic sectional curvature are obviously invariant under the left action of ${\rm U}(p)$ on $T_o^{1,0}\X\simeq M_{p,q}(\C)$. 

Let $\{A_1,\dots, A_d\}$ be an orthonormal basis of $S$. We are going to show that there exists $U\in{\rm U}(p)$ such that

$$UA_k=\frac{1}{2\sqrt q}\left(\begin{array}{c}
0_{(k-1)q,q}\\I_q\\ 0_{p-kq,q}\\
\end{array}\right)\ ,\ k=1,\dots, d.$$
       
Since the column vectors of each $A_i$ are pairwise orthogonal and have the same norm, there exists $U_1\in{\rm U}(p)$ such that

$$U_1A_1=\frac{1}{2\sqrt q}\left(\begin{array}{c}
I_q\\ 0_{p-q,q}\\
\end{array}\right).$$
Now, for any $(\lambda,\mu)\in\C^2\backslash\{(0,0)\}$,
$\trans\hskip 2pt\overline{(\lambda A_1+\mu A_2)}(\lambda A_1+\mu A_2)$ must be a (nonzero) multiple of $I_q$. This implies that for any $(\lambda,\mu)$, $\bar\lambda\mu\,\trans\bar A_1A_2+\lambda\bar\mu\, A_1\hskip 2pt\trans \bar A_2$ is a multiple of $I_q$ (because $\trans\bar A_1A_1$ and $\trans\bar A_2A_2$ are). Moreover, it is trace free because $A_1$ and $A_2$ are orthogonal.
Thus, $\trans\bar A_1A_2=0$ (that is each column vector of $A_1$ is orthogonal to every column vector of $A_2$) and there exists $U_2\in{\rm U}(p)$ such that $U_2U_1A_1=U_1 A_1$ and
$$U_2U_1A_2=\frac{1}{2\sqrt q}\left(\begin{array}{c}
0_{q,q}\\I_q\\ 0_{p-2q,q}\\
\end{array}\right).$$
 
One might continue this process and after $d$ steps, one obtains $U=U_d\dots U_1$. It is then clear that $d$ must be less than or equal to  $p/q$.
\end{demo}

The embedding ($\star$) will be denoted by $f_{max}$ and called the maximal embedding of $\Bm$ into $\X$. This is because ${f_{max}}^\star g_\X=q\,g$, that is, for any $x\in \Bm$ and any $X\in T_x\Bm$, $g_\X({\rm d}f_{max}(X),{\rm d}f_{max}(X))= q\,g(X,X)$, while for a general holomorphic map $f:\Bm\fd \X$ we only know that $f^\star g_\X\leq q\,g$ from the Ahlfors-Schwarz-Pick lemma (see for example~\cite[Theorem~2]{Ro80}). 
Moreover,

\begin{prop}\label{maximal}
Let $f:\Bm\fd \X=\SU(p,q)/{\rm S}(\U(p)\times\U(q))$ be a holomorphic map such that $f^\star g_\X= q\,g$ holds everywhere. Then $p/q\geq m$ and $f$ is totally geodesic. In fact, up to composition with an isometry of $\X$, $f$ is the maximal embedding $f_{max}$.
\end{prop}

\begin{demo}
If $X$ is a nonzero tangent vector at $x\in\Bm$, we denote by $C_{x,X}$ the complex geodesic through $x$ that is tangent to $X$. Let $z$ be a (global) complex coordinate on $C_{x,X}$ and let $\rho^2\,|dz|^2$ (resp. $\sigma^2\,|dz|^2$) be the Hermitian metric induced by $g$ (resp. $f^\star g_\X=q\, g$) on $C_{x,X}$. The Gaussian curvature of $\rho$ (resp. $\sigma$) is given by $K=-\frac{1}{\rho^2}\Delta\log\rho$ (resp. $k=-\frac{1}{\sigma^2}\Delta\log\sigma$). Since $C_{x,X}\subset\Bm$ is totally geodesic, we have $K\equiv -1$. Moreover, because of the holomorphicity of $f$, $k$ is bounded from above by $-1/q$ which is the maximum of the holomorphic sectional curvature on $\X$ and $k\equiv-1/q$ iff the restriction of $f$ to $C_{x,X}$ is totally geodesic. But $\sigma=\sqrt q\,\rho$ and so
$$
k=-\frac{1}{\sigma^2}\,\Delta\log\sigma=-\frac{1}{q\rho^2}\,\Delta\log\rho=\frac{1}{q\rho^2}\,\rho^2K=-\frac{1}{q}.
$$
Thus, $k\equiv-1/q$ and, since this is true for any $(x,X)$, $f$ must be totally geodesic.

Let $o=0_{m,1}$ (resp. $o'=0_{p,q}$) be fixed origins in $\Bm$ (resp. in $\X$). One may suppose (after composition with an isometry of $\X$) that $f(o)=o'$. A consequence of the preceding discussion is that $\d f(T_o\Bm)$ is a $m$-dimensional complex subspace of $T_{o'}\X$ on which the restriction of the holomorphic sectional curvature is constant, equal to $-1/q$. By the proof of Lemma~\ref{tangentmax}, we know that $m\geq p/q$ and that, after composition of $f$ with a suitable isometry, $\d f_{|T_o\Bm}=\d {f_{max}}_{|T_o\Bm}$. By uniqueness of the totally geodesic map satisfying this condition, one has $f=f_{max}$.
\end{demo}

Another maximal embedding of $\Bm$ into $\X=\SU(p,q)/{\rm S}(\U(p)\times\U(q))$ ($p/q\geq m$) is given by
$$
f'_{max}:z=\left(\begin{array}{c}z_1\\ z_2\\ \vdots \\ z_m
        \end{array}\right)\longmapsto
Z=\left(\begin{array}{cccc}z &0 & \cdots& 0\\ 0&z&\cdots&0\\ \vdots & \vdots & \ddots & \vdots\\
0&0&\cdots&z\\ 0&0&\cdots&0\\ \vdots&\vdots&\vdots&\vdots\\ 0&0&\cdots&0\\
        \end{array}\right).
$$
From the previous proposition (it can also be easily verified by hand), this embedding is equal to
$f_{max}$ composed with an isometry of $\X$. The geometric picture is maybe clearer here: $f'_{max}$ is a diagonal embedding of $\B^m$ into $({\B^m})^q\subset\X$ corresponding to a (diagonal) embedding of $\SU(m,1)$ into $\SU(m,1)^q\subset\SU(p,q)$. Moreover the stabilizer in $\SU(p,q)$ of the image of $f'_{max}$
can be computed quite easily. First, let us consider the subgroup of $\U(p,q)$ consisting of
elements of the form:
$$
\left(\begin{array}{ccccccc}
A_{11}&\cdots&A_{1q}& 0& B_{11}&\cdots&B_{1q}\\
\vdots&\cdots&\vdots& 0 & \vdots&\cdots&\vdots\\
A_{q1}&\cdots&A_{qq}& 0& B_{q1}&\cdots&B_{qq}\\
0 &\cdots&0&U&0&\cdots&0\\
C_{11}&\cdots&C_{1q}& 0& d_{11}&\cdots&d_{1q}\\
\vdots&\cdots&\vdots& 0 & \vdots&\cdots&\vdots\\
C_{q1}&\cdots&C_{qq}& 0& d_{q1}&\cdots&d_{qq}\\
\end{array}\right)
$$
where $A_{ij}\in M_{m}(\C)$, $B_{ij}\in M_{m,1}(\C)$, $C_{ij}\in M_{1,m}(\C)$, $d_{ij}\in\C$ and $U\in
\U(r)$ ($r=p-qm$).

Let us denote by $S_q$ the symmetric group on $q$ letters, and by $\U(1)^q\rtimes S_q$ the semi-direct product of 
$\U(1)^q$ by $S_q$ given by the group operation
$(\alpha,\sigma).(\beta,\tau)=(\alpha.\sigma(\beta),\tau\circ\sigma)$.

Define a group homomorphism $\varphi$ of  $\left(\U(1)^q\rtimes S_q\right)\times\SU(m,1)\times\U(r)$ in the above subgroup of $\U(p,q)$ in the following way: if $\alpha=(\alpha_1,\dots,\alpha_q)\in\U(1)^q$, $\sigma\in S_q$, $u\in U(r)$ and
$$
g=\left(\begin{array}{cc}
A & B\\
C & d\\
\end{array}\right)\in \SU(m,1)\ \ \hbox{(where $A\in M_{m}(\C)$, $B\in M_{m,1}(\C)$, $C\in M_{1,m}(\C)$,
$d\in\C$),}
$$
the image of $(\alpha,\sigma,g,u)$ in $\U(p,q)$ is the matrix defined by
$A_{i\sigma(i)}=\alpha_i A$, $B_{i\sigma(i)}=\alpha_i B$, $C_{i\sigma(i)}=\alpha_i C$,
$d_{i\sigma(i)}=\alpha_i d$, $U=u$, and the other coefficients are zero. 

Then ${\rm Ker}\varphi$ is isomorphic to $\Z/(m+1)\Z$ and the stabilizer in $\SU(p,q)$ of the image of $f'_{max}$ is 
${\rm Im}\varphi\cap\SU(p,q)$.

\subsection{Toledo invariant and $\SU(p,2)$-Higgs bundles}\label{sup2-higgs}\hfill

Now we consider a reductive representation $\rho$ of a torsion-free uniform lattice of $\SU(m,1)$, $m>1$, into the Lie group of Hermitian type $G=\SU(p,q)$, $p\geq q\geq 1$. 
Let $M$ be the closed complex hyperbolic manifold $\G\backslash\Bm$. As explained in the introduction, the Toledo invariant can be expressed using the degree of the pull-back of the holomorphic tangent bundle $T^{1,0}\X$ of $\X=G/K$ by any $\rho$-equivariant map $f:\Bm\fd\X$, which we can choose to be harmonic.  
Let then $(P_{K_\C},\t)$ be the $G$-principal Higgs bundle on $M$ associated to $\rho$ and $f$ as in section~\ref{higgs} and  
let $E$ be the holomorphic vector bundle on $M$ associated to $P_{K_\C}$ via the action of $K_\C$ on $\E=\V\oplus\W$. Since $K_\C$ respects the decomposition  $\E=\V\oplus\W$, the bundle $E$ splits holomorphically as the sum of the rank $p$ subbundle $V=P_{K_\C}\times_{K_\C}\V$ with the rank $q$ subbundle $W=P_{K_\C}\times_{K_\C}\W$. 
As a differentiable bundle $E$ is the bundle associated to the flat principal bundle $P_G$ via the action of $G$ on $\E$: it is flat, hence of degree 0. In particular, $\deg V=-\deg W$.

The Higgs field $\t$ is a holomorphic (1,0)-form taking values in the bundle $P_{K_\C}\times_{\Ad K_\C}\pp_\C=\Hom(W,V)\oplus\Hom(V,W)$ so that we can write (see also~\cite{X,BGPG1})
$$
\t=\left(\begin{array}{cc} 0 & \b\\ \g & 0\end{array}\right),\mbox{ where }\left\{\begin{array}{l}\b:W\fd V\otimes\Om\\ \g:V\fd W\otimes\Om \end{array}\right.
$$ 
Recall that the Higgs vector bundle $(E,\t)$ is polystable.  

It is clear that the bundle $f^\star T^{1,0}\X$ is nothing but the bundle $\Hom(W,V)$ and therefore its degree is simply given by $p\,\deg W^\star+q\,\deg V=-(p+q)\deg W$. We obtain that $\tau(\rho)=\frac{4\pi}{m!}\deg W$, so that the Milnor-Wood type inequality that should hold reads
$$
|\deg W| \leq \frac{q}{m+1}\,\deg \Om.
$$
Therefore, for $q=2$ and $\rho$ reductive, our main theorem can be reformulated: 

\begin{theo}\label{ineq} 
Let $\G$ be a torsion free uniform lattice in $\SU(m,1)$, $m>1$, and let $\rho:\G\fd \SU(p,2)$ be a reductive representation. Let $E=V\oplus W$ be the $\SU(p,2)$-Higgs vector bundle on $M=\G\backslash\Bm$ associated to $\rho$. Then 
$|\deg W| \leq \frac{2}{m+1}\,\deg \Om$ with equality if and only if  $m\leq p/2$ and, up to conjugacy, $\rho$ is induced by the maximal embedding $f_{max}:\Bm\fd\X$ or by its conjugate. 
\end{theo}

If we deform the Higgs bundle $(E,\t)$ via the $\C^\star$-action on the moduli space as in section~\ref{higgs}, we obtain a system of Hodge bundles:  
$$
E_1\stackrel{\t_1}{\fd}E_2\stackrel{\t_2}{\fd}\ldots\stackrel{\t_{k-1}}{\fd}E_k\stackrel{\t_k}{\fd} 0.
$$ 
Moreover, each subbundle $E_i$ splits as $E_i=V_i\oplus W_i$ with $V_i\subset V$ and $W_i\subset W$, and $\t_i$ decomposes as $\g_i\oplus\b_i$, where $\g_i:V_i\fd W_{i+1}\otimes\Om$ and $\b_i:W_i\fd V_{i+1}\otimes\Om$. We obtain two Higgs subbundles 
$$
V_1\fd W_2 \fd V_3 \fd W_4 \fd \ldots \fd 0
$$ 
and 
$$
W_1\fd V_2 \fd W_3 \fd V_4 \fd \ldots \fd 0
$$
which are again polystable of degree 0.   

So far everything we said was valid in the general rank $q$ case. Now we will need the assumption that $q=2$ to ensure that the systems of Hodge bundles we obtain are simple ones and/or that the decomposition of $W$ is maximal in the sense that $W$ splits into a sum of line bundles.  
Indeed, if $W$ has rank 2, we see that by regrouping and renaming the subbundles if necessary, we can write our system of Hodge bundles either as a polystable Higgs bundle of the form 
$$
V_1\stackrel{\g_1}{\fd}W\stackrel{\b}{\fd}V_2\stackrel{\g_2}{\fd}0 
$$
with $V_1\oplus V_2=V$, or as a polystable Higgs bundle of the form
$$
V_1\stackrel{\g_1}{\fd}W_1\stackrel{\b_1}{\fd}V_2\stackrel{\g_2}{\fd}W_2\stackrel{\b_2}{\fd}V_3\stackrel{\g_3}{\fd}0
$$
where $W_1$ and $W_2$ are line bundles, $W_1\oplus W_2=W$ and $V_1\oplus V_2\oplus V_3=V$. 

\subsection{Proof of the Milnor-Wood type inequality}\label{inequality}\hfill

The case of non reductive representations will be postponed to Paragraph~\ref{nonreductive}. Until then, the representation $\rho:G\fd \SU(p,2)$ is assumed to be reductive so that we can apply the results of Section~\ref{sup2-higgs}. Our proof of the inequality $|\deg W|\leq\frac{2}{m+1}\,\deg\Om$ will be different according to the form of the system of Hodge bundles we obtain by deforming the polystable Higgs bundles $E=V\oplus W$ via the $\C^\star$-action. Note that the deformation changes the holomorphic structures of $E$, $V$ and $W$, but not their isomorphism classes as differentiable complex vector bundles, hence their degrees remain unchanged. During the proof, we will see that if equality holds, some bundles have stability properties (in the usual sense) that will be useful for the study of maximal representations in Section~\ref{equality}. 

We refer to~\cite{VZ05,HL97} for general facts about sheaves and stability.

\subsubsection{System of Hodge bundles of the type $V_1{\fd}W{\fd}V_2{\fd}0$}\label{1W}\hfill

Here the important point is that the system of Hodge bundles we are dealing with is a ternary bundle, and no limitation on the rank of $W$ is needed. Hence the results of this paragraph are valid in the general case $\rk W=q\geq 1$. 

Let ${\cal F}$ be the maximal destabilizing subsheaf of $W$, that is, the first term in the Harder-Narasimhan filtration of $W$~\cite{VZ05,HL97}. By definition, ${\F}$ has maximal slope among the subsheaves of $W$, hence is semistable.   Consider the restriction $\b_{\cal F}:{\cal F}\otimes\T\fd V_2$. Since $\t$ vanishes on $V_2$, ${\cal F}\oplus\Im\b_{\cal F}$ is a Higgs subsheaf and hence by stability, $\deg \Im\b_{\cal F}\leq -\deg {\cal F}$. Now, the tensor product of two semistable sheaves is again semistable and hence ${\cal F}\otimes\T$ is semistable. Therefore we have  $\mu({\cal F})+\mu(\T)=\mu({\cal F}\otimes\T)\leq\mu(\Im\b_{\cal F})$ which implies $(\rk\b_{\cal F}+\rk {\cal F})\mu({\cal F})\leq\rk\b_{\cal F}\mu(\Om)$. Thus,
$$
\deg W\leq q\mu({\cal F})\leq q\frac{\rk \b_{\cal F}}{\rk \b_{\cal F} + \rk {\cal F}}\frac{\deg \Om}{m}\leq\frac{q}{m+1}\deg\Om
$$
where the first inequality follows from the fact that ${\cal F}$ is of maximal slope among the subsheaves in $W$, and the last from $\rk \b_{\cal F}\leq m\rk {\cal F}$.  

The remaining inequality is obtained exactly in the same way by considering the dual Higgs bundle 
$$
V_2^\star\stackrel{\trans\b}{\fd}W^\star\stackrel{\trans\g_1}{\fd}V_1^\star{\fd}0. 
$$

Assume that equality holds, for example that $\deg W=\frac{q}{m+1}\deg\Om$. Then, retracing our steps, we see that $W$ must be a semistable bundle (in the usual sense), that we must have $\rk \b=mq$, i.e. $\b:W\otimes\T\fd V$ generically injective, and moreover that $\deg W\oplus\Im\b=0$. This last fact implies by polystability that $E$ splits as the sum of $(W\oplus\Im \b,\b)$ with an other polystable Higgs bundle $E'$ of degree 0. In our case, this means that $V_2$ splits holomorphically as $\Im \b\oplus V_2'$ and that $\g_1$ vanishes. $E'$ is then the polystable (in the usual sense) subbundle $V_1\oplus V_2'$ of $V$. 

In the same manner,  we find that $\deg V=\frac{q}{m+1}\deg\Om$ implies that $W$ is a semistable bundle, $\b=0$, $V_1$ splits holomorphically as $\Ker \g_1\oplus V_1'$, and $\g_1:V_1'\fd W\otimes\Om$ is an isomorphism. Our system of Hodge bundle is the sum of the Higgs bundle $(V_1'\oplus W,\g_1)$ with the polystable subbundle $\Ker\g_1\oplus V_2$ of $V$.   

\subsubsection{System of Hodge bundles of the type $V_1{\fd}W_1{\fd}V_2{\fd}W_2\fd V_3\fd 0$}\label{2W}\hfill

Here, we need to assume that $\rk W=q=2$, namely that $W_1$ and $W_2$ are line bundles. 

Assume first that $\b_1:W_1\otimes\T\fd V_{2}$ vanishes. 
We then have to deal with the sum of two polystable Higgs bundles of degree 0: $V_1\fd W_1\fd 0$ and $V_2\fd W_2 \fd V_3\fd 0$. In this situation, we already know from \ref{1W} that $0\leq\deg V_1\leq\frac{1}{m+1}\deg\Om$ and $|\deg (V_2\oplus V_3)|\leq \frac{1}{m+1}\deg\Om$ which gives the result.

We are left with the case where  $\b_1:W_1\otimes\T\fd V_{2}$ is non zero. 

Consider  $\g_1:V_1\fd W_1\otimes\Om$. We have $\deg V_1=\deg\Ker\g_1+\deg\Im\g_1\leq \rk\g_1\,\mu (W_1\otimes\Om)$ since $\Ker\g_1$ is $\t$-invariant 
and $W_1\otimes\Om$ is a semistable bundle, being the product of a stable bundle by a line bundle. Since $\deg W_1=-(\deg V_1+\deg V_2+\deg W_2+\deg V_3)$ we obtain
$$
\deg V_1\leq \frac{\rk\g_1}{1+\rk\g_1}\frac{\deg\Om}{m}-\frac{\rk\g_1}{1+\rk\g_1}(\deg V_2+\deg W_2+\deg V_3).
$$
In the same way, 
$$
\deg V_2\leq \frac{\rk\g_2}{1+\rk\g_2}\frac{\deg\Om}{m}-\frac{\rk\g_2}{1+\rk\g_2}(\deg V_1+\deg W_1+\deg V_3).
$$
Hence,
$$
\begin{array}{rcl}
\deg V & \leq & \ds \left(\frac{\rk\g_1}{1+\rk\g_1}+\frac{\rk\g_2}{1+\rk\g_2}\right)\frac{\deg\Om}{m}\, + \\
&& \ds \left(\frac{\rk\g_2}{1+\rk\g_2}-\frac{\rk\g_1}{1+\rk\g_1}\right)(\deg V_2+\deg W_2+\deg V_3)
+\left(1-\frac{\rk\g_2}{1+\rk\g_2}\right)\deg V_3.
\end{array}
$$ 

Now the commutation relation $[\t,\t]=0$ gives us control on the rank of the $\g_i$'s: 

\begin{lemma}
Assume that $\b_i:W_i\otimes\T\fd V_{i+1}$ is non zero. Then the rank of $\g_i:V_i\fd W_i\otimes\Om$ is at most 1. 
\end{lemma}

\begin{demo}
This is linear algebra. We work in a single fiber above some point in $M$. We write $W_i=\C w_i$. Since $\b_i$ is non zero, there exists $Z\in\T$ such that $\b_i(Z)w_i\neq 0$. Assume that the rank of $\g_i$ is at least 2 at some point. Then we can find two linearly independent forms $\a,\a'$ and two vectors $v,v'$ in $V_i$ such that 
$\g_i(v)=w_i\otimes\a$ and $\g_i(v')=w_i\otimes\a'$. Now, the commutation relation $[\t,\t]=0$ means in particular that for all $X,Y\in\T$ and all $u\in V$ , $\b_i(X)\g_i(Y)u=\b_i(Y)\g_i(X)u$. 

If we take $u=v$, $Y=Z$ and $X$ such that $\a(X)\neq 0$, we get $\a(Z)\b_i(X)w_i=\a(X)\b_i(Z)w_i$, which implies that $\a(Z)$ and 
$\b_i(X)w_i$ are different from zero. So if now $X'$ is such that $\a(X')=0$, $\b_i(X')w_i=0$. We can choose such an $X'$ with the additional property that $\a'(X')\neq 0$, since $\a$ and $\a'$ are independent. This is a contradiction since we could have taken $u=v'$ to prove that $\b_i(X')w_i\neq 0$ if $\a'(X')\neq 0$.    
\end{demo}

We assumed that $\b_1:W_1\otimes\T\fd V_{2}$ is non zero, hence $\rk\g_1\leq 1$. If $\rk\g_2\geq \rk\g_1$, using the fact that
$\deg V_2+\deg W_2+\deg V_3$ and $\deg V_3$ are both non positive, and since $\rk\g_2\leq m$, we get 
$$
\deg V\leq \frac{3m+1}{2m(m+1)}\deg\Om\leq\frac{2}{m+1}\deg\Om,
$$ 
and the last inequality is strict as soon as $m>1$. 
 If $\rk\g_2 < \rk\g_1$, that is if $\rk\g_1=1$ and  $\rk\g_2=0$, we are in the case $V_1\fd W_1\fd V_2\fd 0$ and hence 
 $$\deg V\leq \deg(V_1\oplus V_2)=-\deg W_1\leq\frac{1}{m+1}\deg\Om.$$  

Once again, the remaining inequalities are obtained by looking at the dual Higgs bundle 
$$
V_3^\star\stackrel{\trans\b_2}{\fd}W_2^\star\stackrel{\trans\g_2}{\fd}V_2^\star\stackrel{\trans\b_1}{\fd}W_1^\star\stackrel{\trans\g_1}{\fd}V_1^\star{\fd}0. 
$$ 

\vspace{0.3cm}

Assume we are in the equality case and $m>1$. If $\deg V= \frac{2}{m+1}\deg\Om$, $\b_1:W_1\otimes\T\fd V_2$ vanishes, and the equalities $\deg V_1=\frac{1}{m+1}\deg\Om=\deg(V_2\oplus V_3)$ hold. We saw in \ref{1W} that in this situation, $\b_2=0$ and there exists holomorphic subbundles $V_i'\subset V_i$ such that $\g_i:V_i'\fd W_i\otimes \Om$ are isomorphisms for $i=1,2$.  

If $\deg W=\frac{2}{m+1}\deg\Om$, we find that for $i=1,2$, $\g_i=0$ and there exists a holomorphic subbundle $V_i'\subset V_i$ such that $\b_i:W_i\otimes\T\fd V_i'$ is an isomorphism.

In either cases, $\deg W_1=\deg W_2$ and $W$ with the deformed complex structure is polystable hence semistable. 

\subsubsection{Non reductive representations}\label{nonreductive}\hfill

Assume now that the representation $\rho:\G\fd G=\SU(p,2)$ is not reductive. This implies that $\rho(\G)$ fixes a point $\xi$ on the boundary at infinity $\X(\infty)$ of $\X$ (\cite{Labourie}). Let us fix an origin $o\in\X$ and let $c$ be the unit speed geodesic ray starting from $o$ representing $\xi$. Let $\gg=\kk\oplus\pp$ be the Cartan decomposition of $\gg$ associated to $o$ and let $X\in\pp$ be such that $X=\dot c(0)$ in the usual identification of $\pp$ with $T_o\X$. We have the following description of the stabilizer $G_\xi$ of $\xi$ in $G$ (see for example~\cite{Eberlein}):
$$
G_\xi=\{g\in G\,|\,\lim_{t\rightarrow +\infty}\exp(-tX)\,g\,\exp(tX)\mbox{ exists}\}=K_\xi.A_\xi.N_\xi
$$   
where $A_\xi=\exp(\{Y\in\pp\,|\,[X,Y]=0\})$, $N_\xi= \displaystyle\{g\in G_\xi\,|\,\lim_{t\rightarrow +\infty}\exp(-tX)\,g\,\exp(tX)=1\}$, and  $K_\xi=G_\xi\cap K$.  

By assumption $\rho(\G)\subset G_\xi$ and we can consider the so-called semi-simplification $\rho_{ss}$ of $\rho$ which is defined by $\rho_{ss}(\g)= \lim_{t\rightarrow +\infty}\exp(-tX)\,\rho(\g)\,\exp(tX)\in K_\xi.A_\xi$ for all $\g\in\G$. The representation $\rho_{ss}$ belongs to the connected component of $\rho$ in the space $\Hom(\G,G)$ and is reductive: we can apply the results of the last paragraphs to get the Milnor-Wood bound on $\tau(\rho)=\tau(\rho_{ss})$. 

In fact we can do better. The representation $\rho_{ss}$ stabilizes the orbit $K_\xi.A_\xi.o=A_\xi.o$, which is a totally geodesic submanifold of $\X$. It is not difficult to see (and probably well known) that this orbit is either a totally real totally geodesic submanifold of $\X$ (for example, if $c$ is a regular geodesic, it is the unique maximal flat, isometric to $\R^2$ in our case, containing $c$), or the Riemannian product of $\R$ with a totally geodesic copy of complex hyperbolic space $\B^{p-1}$ (of induced holomorphic sectional curvature $-1$). In the first case the Toledo invariant is zero since the restriction of the Kähler form $\o_\X$ to a totally real submanifold vanishes. In the second one it is bounded (in absolute value) by ${\rm Vol}(M)$.     

Therefore non reductive representations can not be maximal. 

\subsection{Maximal representations}\label{equality}\hfill

Thanks to the previous paragraph, we know that if the representation $\rho$ is maximal, it is reductive. Therefore we may consider the polystable Higgs bundles $(E=V\oplus W,\t)$  associated to $\rho$. 

In order to prove the Milnor-Wood type inequality $|\deg W|\leq \frac{2}{m+1}\,\deg\Om$, we have deformed the Higgs bundle $E$ to a system of Hodge bundles. 
Here, we need to distinguish between these two Higgs bundles, and we will call 
the latter $(E_0=V_0\oplus W_0,\t_0)$. Let $\dbar_{W}$  and $\dbar_{W_0}$ be the complex structure of $W$ and $W_0$. Again, although the complex structure is (by definition) not modified by the $\C^\star$-action, in the limit $\dbar_{W_0}$ is a priori different from $\dbar_W$. In fact, all we know is that $\dbar_{W_0}$ is in the closure of the orbit of $\dbar_W$ under the group of gauge transformations: i.e. there exist gauge transformations $g_{t_i}$ such that $g_{t_i}^\star\dbar_W$ goes to $\dbar_{W_0}$ when $t_i$ goes to 0. Let us call $W_{t_i}$ the bundle $W$ with the complex structure $g_{t_i}^\star\dbar_W$. 

The main point of Section~\ref{inequality}, apart from the proof of the inequality itself, was that $|\deg W|=\frac{2}{m+1}\,\deg\Om$ implies that the bundle $W_0$ is semistable (in the usual sense). This implies that $W$ itself is semistable (regardless of the rank of $W$):

\begin{lemma}
Assume that $W$ with its initial complex structure $\dbar_W$ is not a semistable bundle. Then $W_0$, that is $W$ with the complex structure $\dbar_{W_0}$, is not semistable either.  
\end{lemma}

\begin{demo}
Let ${\cal F}$ be a subsheaf of $(W,\dbar_W)$ such that $\mu({\cal F})>\mu(W)$. Let $r$ be the rank of ${\cal F}$. We have a monomorphism of sheaves ${\cal F}\fd W$ and therefore, for all $i\in\N$, we obtain a monomorphism  of sheaves ${\cal F}\fd W_{t_i}$. This gives non trivial holomorphic maps between the determinant bundle $\det{\cal F}=(\bigwedge^r{\cal F})^{\star\star}$ of ${\cal F}$ and $\bigwedge^r W_{t_i}$. This means that for all $i\in\N$, the cohomology group $H^0(M,\Hom(\det{\cal F},\bigwedge^r W_{t_i}))$ is at least one dimensional. By the upper semicontinuity of cohomology (see Kobayashi~\cite{K}), there exists a non-trivial holomorphic map  $\det{\cal F}\fd\bigwedge^r W_0$. Let $I$ be its image. Since $\det{\cal F}$ is stable (it is a line bundle), we have 
$$
\mu(I)\geq\mu(\det{\cal F})=r\mu({\cal F})> r\mu(W)=r\mu(W_0)=\mu(\bigwedge ^r W_0).
$$ 
Hence $\bigwedge^r W_0$ is not semistable, and neither is $W_0$.    
\end{demo}

Summing up, we proved 

\begin{prop}
Let $\rho:\G\fd\SU(p,2)$ be a maximal representation and let $E=V\oplus W$ be the 
associated $\SU(p,2)$-Higgs bundle on $M=\G\backslash\Bm$. Then $W$ is a semistable holomorphic bundle.
\end{prop}

The semistability of $W$ is a very strong property and Theorem~\ref{ineq} follows from 

\begin{theo}
Let $\G$ be a torsion-free uniform lattice in $\SU(m,1)$, $m>1$, and $\rho:\G\fd\SU(p,q)$, $p\geq q\geq 1$, be a reductive representation.  Let $E=V\oplus W$ be the associated $\SU(p,q)$-Higgs bundle on $M=\G\backslash\Bm$. Assume moreover that $W$ is semistable. Then $|\deg W|\leq\frac{q}{m+1}\,\deg\Om$,  with equality if and only if $m\leq p/q$ and, up to conjugacy, $\rho$ is induced by the maximal embedding $f_{max}:\Bm\fd\X$ or by its conjugate.
\end{theo}
 
\begin{demo}
Consider $\b:W\otimes\T\fd V$ and argue as in ~\ref{1W} with $\F=W$ and $V_2=V$ to get the bound $\deg W\leq\frac{q}{m+1}\,\deg\Om$, with equality if $\b$ is injective and $W\oplus\Im\b$ has degree zero and hence is a polystable Higgs subbundle of $E$. Now, we have the 

\begin{lemma}
If $m>1$ and $\b:W\otimes\T\fd V$ is injective, $\g$ vanishes identically. 
\end{lemma}

\begin{demo}[of the lemma]
This is again a consequence of the relation $[\t,\t]=0$, which in our case reads $\b(X)\g(Y)v=\b(Y)\g(X)v$ for all $X,Y\in\T$ and all $v\in V$. Let $\{w_1\ldots,w_q\}$ be a basis of $W$ above some point $x\in M$. We can write 
$\g(X)v=\sum_{i=1}^q \lambda_i(X,v)w_i$ and $\g(Y)v=\sum_{i=1}^q \lambda_i(Y,v)w_i$. But this implies that 
$$
\sum_{i=1}^q \lambda_i(Y,v)\b(X)w_i=\sum_{i=1}^q \lambda_i(X,v)\b(Y)w_i.
$$
Since $m>1$, we can take $X$ and $Y$ to be linearly independent and the injectivity of $\b$ implies that $\b(X)w_1,\ldots,\b(X)w_q,\b(Y)w_1,\ldots,\b(Y)w_q$ are linearly independent vectors in $V$, hence that $\lambda_i(X,v)=\lambda_i(Y,v)=0$ for all $i$, namely that $\g=0$.  \end{demo}

Therefore $\deg W=\frac{q}{m+1}\,\deg\Om$ implies that $\g=0$, hence that $\partial^{0,1}f=0$, i.e. the harmonic map $f$ is holomorphic. The theorem easily follows. We know from the Alhfors-Schwarz-Pick lemma that $f^\star g_\X\leq q\,g$. But this implies that the inequality $\la f^\star\o_\X,\o\ra\leq 2mq$ is pointwise true whereas $\deg W=\frac{q}{m+1}\deg\Om$ means that  $\tau(\rho)=\frac{1}{2m}\int_M\la f^\star\o_\X,\o\ra\,dV = q{\rm Vol}(M)$, so that in fact $f^\star g_\X= q\,g$ holds everywhere. Proposition~\ref{maximal} yields that $f=f_{max}$, up to composition with an isometry of $\X$.

To get the inequality $\deg W\geq -\frac{q}{m+1}\,\deg\Om$, consider the map  
$\g:V\fd W\otimes\Om$. We have $\deg V=\deg\Ker\g+\deg\Im\g$. Since $\Ker\g$ is $\t$-invariant, $\deg\Ker\g\leq 0$. By semistability of $W$, $\deg\Im\g\leq \rk\g(\frac{1}{q}\deg W+\frac{1}{m}\deg \Om)$. Hence 
$$
\deg V\leq\frac{q\,\rk \g}{q+\rk\g}\frac{\deg\Om}{m}\leq\frac{q}{m+1}\deg\Om
$$
with equality if and only if $\rk\g=qm$, i.e. $\g$ is generically onto, and $\deg\Ker\g=0$, i.e. $\Ker\g$ is a polystable subbundle of $E$. 

Again, the fact that $[\t,\t]=0$ yields that $\b=0$, i.e. $f$ is antiholomorphic:

\begin{lemma}
If $m>1$ and $\g:V\fd W\otimes\Om$ is onto, $\b$ vanishes identically.
\end{lemma}  
   
\begin{demo}[of the lemma]
Let $X\in\T$ and $w\in W$. Let $\alpha\in\Om$, $\a\neq 0$, be such that $\alpha(X)= 0$. Take $v\in V$ such that $\g(v)=w\otimes\alpha$. Then for all $Y\in \T$, we have on the one hand $\b(X)\g(Y)v=\b(X)(\alpha(Y)w)=\a(Y)\b(X)w$ and on the other hand $\b(X)\g(Y)v=\b(Y)\g(X)v=\a(X)\b(Y)w=0$. We may find $Y$ such that $\alpha(Y)\neq 0$, for $m$ is greater than 1. Hence $\b(X)w=0$.  
\end{demo}

The rest of the proof goes like in the holomorphic case.
\end{demo}

\subsection{Proof of Proposition~\ref{otherbound}}\hfill

We use freely what has been done in Sections~\ref{sup2-higgs} and~\ref{inequality}. If the representation $\rho$ is not reductive, we consider its semi-simplification instead. Considering the polystable Higgs bundle $E=V\oplus W$ associated to $\rho$, we want to prove the inequality
$$
|\deg W|\leq\frac{2p}{p+2}\,\frac{\deg\Om}{m}.
$$  
The proof again depends on the type of system of Hodge bundle we obtain by deforming $E$ via the $\C^\star$-action. 

\subsubsection{System of Hodge bundles of type $V_1\fd W\fd V_2\fd 0$}\hfill 

As in Paragraph~\ref{1W}, we need no restriction on the rank of $W$ here. So let $q=\rk W\geq 1$.

We use the method of Viehweg and Zuo~\cite{VZ05}. 
They work with a binary system of Hodge bundles ($V_1=0$) so we explain how to adapt their proof to the ternary case. We try to fit to their notations as much as possible. Dualizing the Higgs bundle if necessary, we may suppose that $\deg W>0$. We also suppose that no subsheaf of $V_2$ has a slope equal to zero. In fact, each subsheaf of $V_2$ has non positive slope because $\theta_{|{V_2}}=0$ and if its slope is equal to zero, then the Higgs bundle splits as a sum of two polystable Higgs bundles of degree zero with one contained in $V_2$.

Let us consider the Harder-Narasimhan filtrations~\cite{VZ05,HL97} 
$$0=W^0\subset W^1\subset\dots\subset W^{l''}=W$$
and
$$0=V^0_2\subset V^1_2\subset\dots\subset V^{l'}_2=V_2$$
of $W$ and $V_2$. Let $l$ be the maximum of all $j$ verifying $\mu(W^j/W^{j-1})>0$. 
Remark that $l\geq 1$ because $\mu(W^1)\geq\mu(W)>0$.

We construct by induction two sequences
$$0=j_0<j_1<\dots<j_r=l\ \ {\rm and}\ \ 0=j'_0<j'_1<\dots<j'_r\leq l'$$
in the following way:

\noindent Suppose that  $j_{k-1}$ and $j'_{k-1}$ are defined. If $j_{k-1}<l$, let $j'_k$ be the minimal number with $\b(W^{j_{k-1}+1})\subset V_2^{j'_k}\otimes\Om$, and $j_k$ be the maximum of all $j\leq l$ verifying $\b(W^j)\subset V_2^{j'_k}\otimes\Om$.


Then, we have non trivial morphisms
$$\frac{W^{j_{k-1}+1}}{W^{j_{k-1}}}\fd\frac{V_2^{j'_k}}{V_2^{j'_k-1}}\otimes\Om.$$
Because of the semistability of all involved sheaves, we get
$$\mu\left(\frac{W^{j_{k-1}+1}}{W^{j_{k-1}}}\right)\leq \mu\left(\frac{V_2^{j'_k}}{V_2^{j'_k-1}}\right)+\mu(\Om)$$
for each $k$. We set $E^k=W^{j_k}\oplus V_2^{j'_k}$. The sequence $(E^k)_{0\leq k\leq r}$ defines a filtration of $W^{j_r}\oplus V_2^{j'_r}$ by Higgs subsheaves, and we denote the successive quotients by $F^k=E^k/E^{k-1}=F_W^k\oplus F_{V_2}^k$, where $F_W^k=W^{j_k}/W^{j_{k-1}}$ and $F_{V_2}^k=V_2^{j'_k}/V_2^{j'_{k-1}}$. 

From the properties of the Harder-Narasimhan filtrations, we have
$$
\mu(F^{k-1}_W)=\mu\left(\frac{W^{j_{k-1}}}{W^{j_{k-2}}}\right)\geq\mu\left(\frac{W^{j_{k-1}}}{W^{j_{k-1}-1}}\right)>\mu\left(\frac{W^{j_{k-1}+1}}{W^{j_{k-1}}}\right)\geq\mu\left(\frac{W^{j_{k}}}{W^{j_{k-1}}}\right)=\mu(F_W^k)
$$
and
$$
\mu(F_{V_2}^k)=\mu\left(\frac{V_2^{j'_k}}{V_2^{j'_{k-1}}}\right)\geq \mu\left(\frac{V_2^{j'_k}}{V_2^{j'_k-1}}\right)>\mu\left(\frac{V_2^{j'_k+1}}{V_2^{j'_k}}\right)\geq\mu\left(\frac{V_2^{j'_{k+1}}}{V_2^{j'_k}}\right)=\mu(F_{V_2}^{k+1}).
$$
In particular, we get for all $1\leq k\leq r$,
$$\mu(F_W^k)-\mu(F_{V_2}^k)\leq \mu(\Om)$$
and
$$\mu(F_W^1)>\mu(F_W^2)>\dots>\mu(F_W^r)>0>\mu(F_{V_2}^1)>\mu(F_{V_2}^2)>\dots>\mu(F_{V_2}^r).$$

Viehweg and Zuo then define the following quantities:

$\bullet$ $c_k=\deg F^k$,

$\bullet$ $\mu_k^W=\mu(F^k_W)$, $\mu_k^{V_2}=\mu(F^k_{V_2})$,

$\bullet$ $r_k^W=\rk(F^k_W)$,

$\bullet$ $r_k^{V_2}=\rk(F^k_{V_2})-\displaystyle\frac{c_k}{\mu_k^{V_2}}$.

With these definitions, we can write the above inequalities
$$\mu_k^W-\mu_k^{V_2}\leq\mu(\Om)\ \ {\rm for\ all\ }1\leq k\leq r,\ \ {\rm and}\ \ \mu_1^W>\dots>\mu_r^W>0>\mu_1^{V_2}>\dots>\mu_r^{V_2}.$$

We verify the properties I--IV of Claim 2.2 in~\cite{VZ05} (recall that in the present situation, we may have $j_r<l''$):

I. Each $r_k^{V_2}=-\frac{\mu_k^W r_k^W}{\mu_k^{V_2}}$ and hence is positive.

II. Here we need some adaptations. For each $0\leq k\leq r$, the number $\sum_{i=1}^k c_i$ is non positive because $E^k$ is a Higgs subsheaf. Moreover (if we set $p_1=\rk V_1$)

$$\displaylines{
\qquad p_1+\sum_{k=1}^r r_k^{V_2}+\rk(V_2/V_2^{j'_r})-p=\sum_{k=1}^r r_k^{V_2}-\sum_{k=1}^r \rk(F^k_{V_2})=-\sum_{k=1}^r \frac{c_k}{\mu_k^{V_2}}\hfill\cr
\hfill{}=-\frac{1}{\mu_r^{V_2}}\Bigl(\sum_{i=1}^r c_i\Bigr)+\sum_{k=1}^{r-1}\frac{\mu_k^{V_2}-\mu_{k+1}^{V_2}}{\mu_k^{V_2}\,\mu_{k+1}^{V_2}}\,\Bigl(\sum_{i=1}^k c_i\Bigr)\leq 0.
\qquad\cr}
$$
Thus, $p\geq p_1+\sum_{k=1}^r r_k^{V_2}+\rk(V_2/V_2^{j'_r})\geq \sum_{k=1}^r r_k^{V_2}$.

III. By assumption, $\deg V<0$ and therefore
$$\begin{array}{rcl}
\mu(V)&\geq&\displaystyle\frac{\sum_{k=1}^r\rk(F^k_{V_2})\mu_k^{V_2}+\deg(V_2/V_2^{j'_r})+\deg V_1}{\sum_{k=1}^r r_k^{V_2}}\\
&=&\displaystyle\frac{\sum_{k=1}^r r_k^{V_2}\mu_k^{V_2}+\sum_{k=1}^r c_k+\deg(V_2/V_2^{j'_r})+\deg V_1}{\sum_{k=1}^r r_k^{V_2}}\\
&=& \displaystyle\frac{\sum_{k=1}^r r_k^{V_2}\mu_k^{V_2}-\deg(W/W^{l})}{\sum_{k=1}^r r_k^{V_2}}\\
&\geq& \displaystyle\frac{\sum_{k=1}^r r_k^{V_2}\mu_k^{V_2}}{\sum_{k=1}^r r_k^{V_2}}
\end{array}
$$
(we use $\sum_{k=1}^r c_k+\deg(V_2/V_2^{j'_r})+\deg V_1=-\deg(W/W^{l})\geq0$).

IV. From III, we get
$$\mu(W)-\mu(V)\leq \frac{\sum_{k=1}^r r_k^{W}\mu_k^{W}}{\sum_{k=1}^r r_k^{W}}-\frac{\sum_{k=1}^r r_k^{V_2}\mu_k^{V_2}}{\sum_{k=1}^r r_k^{V_2}}.$$
The r.h.s is bounded from above by ${\rm max}\bigl\{\mu_k^W-\mu_k^{V_2}\,,\,1\leq k\leq r\bigr\}$ (see~\cite{VZ05} for the proof of this), hence 
$$\frac{p+q}{pq}\,\deg W=\mu(W)-\mu(V)\leq\frac{\deg\Om}{m}.$$

\subsubsection{System of Hodge bundles of the type $V_1{\fd}W_1{\fd}V_2{\fd}W_2\fd V_3\fd 0$}\hfill
 
If $\b_1\neq 0$ and $\g_2\neq 0$, we have exactly as in Section~\ref{2W} that 
$$
|\deg W|\leq\frac{3\rk V_2+1}{2(\rk V_2+1)}\,\frac{\deg\Om}{m}
$$ 
using that $\rk\g_2\leq\rk V_2$ instead of $\rk\g_2\leq m$. This inequality is stronger than the one we want to prove here.

If $\b_1=0$ or $\g_2= 0$, the Higgs bundle splits as the sum of two polystable Higgs bundles of degree 0, for example $V_1\fd W_1\fd 0$ and $V_2\fd W_2\fd V_3\fd 0$. Since $\rk W_1=\rk W_2=1$, we have by the previous paragraph that  
$$
|\deg V_1|\leq \frac{\rk V_1}{\rk V_1+1}\,\frac{\deg\Om}{m} \quad\mbox{ and }\quad |\deg (V_2\oplus V_3)|\leq \frac{\rk V_2+\rk V_3}{\rk V_2+\rk V_3+1}\,\frac{\deg\Om}{m}
$$ 
which, in view of the following lemma, gives the result.

\begin{lemma}
Let $p_1,p_2,q_1,q_2$ be positive numbers. Let $p=p_1+p_2$ and $q=q_1+q_2$. Then
$$\frac{p_1q_1}{p_1+q_1}+\frac{p_2q_2}{p_2+q_2}\leq\frac{pq}{p+q}$$
with equality iff $p_1q_2=p_2q_1$.
\end{lemma}

\begin{demo}
$$
\frac{pq}{p+q}-\frac{p_1q_1}{p_1+q_1}-\frac{p_2q_2}{p_2+q_2}=\frac{(p_1q_2-p_2q_1)^2}{(p+q)(p_1+q_1)(p_2+q_2)}.
$$
\end{demo}

\section{\texorpdfstring{The case $G=\SO_0(p,2)$}{The case G=SO(p,2)}}\label{sop2}

In this section we prove the main theorem in the case where $G=\SO_0(p,2)$, the identity component of $\SO(p,2)$.

\subsection{\texorpdfstring{The Hermitian symmetric space $\SO_0(p,2)/(\SO(p)\times\SO(2))$}{The Hermitian symmetric space SO(p,2)/(SO(p)xSO(2))}}\hfill

Here, the symmetric space $\X$ is naturally seen as an open subset of the (real) Grassmannian of 2-planes of a real vector space, which makes the complex structure a bit more difficult to understand. Again, details are to be found in~\cite{Helgason,Satake}.

Let $\E_\R$ be a real vector space of dimension $p+2$ ($p\geq 3$), endowed with a non-degenerate quadratic form $S$ of signature $(p,2)$. The symmetric space $\X$ is defined as the space of all $2$-dimensional real subspaces of $\mathbb E_\R$ on which $S$ is negative-definite. It is an open submanifold of the real Grassmannian of $2$-planes of $\mathbb E_\R$.

Let us fix a 2-plane $\W_\R$ of $\E_\R$ on which $S$ is negative-definite and let $\V_\R$ be its orthogonal complement. We also fix an orientation on $\W_\R$.

The group $G={\SO}_0(p,2)$ acts transitively on $\X$ by analytic isomorphisms.  The isotropy subgroup $K$ of $G$ at $\W_\R$ is identified with the maximal compact subgroup $\SO(p)\times\SO(2)$, and hence  $\X$ can be identified with ${\SO}_0(p,2)/(\SO(p)\times\SO(2))$.

Let $\gg$ be the Lie algebra of $G$, $\kk\subset \gg$ the Lie algebra of $K$
and $\gg=\kk\oplus \pp$ the corresponding Cartan decomposition. Let us fix an orthonormal basis $(e_1,\dots,e_p)$ of $\V_\R$ and a direct orthonormal basis $(e_{p+1},e_{p+2})$ of $\W_\R$ (with respect to $S_{|\V_\R}$, respectively $S_{|\W_\R}$). We have the following matrix expressions:

$$
\got k=\left\{\left(\begin{array}{cc} X_1 & 0\\ 0 & X_2
\end{array}\right)\ ,\ X_1\in M_p(\R),\ X_2\in M_2(\R),\ \trans X_i=-X_i\ (i=1,2)\right\},
$$ 
$$
\got p=\left\{\left(\begin{array}{cc} 0 & A\\ \trans A & 0
\end{array}\right)\ ,\ A\in M_{p,2}(\R)\right\}\,\simeq_\R\,\Hom(\W_\R,\V_\R).
$$ 

The tangent space $T_o\X$ at $o=\W_\R\in\X$ will be identified with $\got p$. The tangent bundle $T\X$ of $\X$ is the bundle $G\times_{\Ad K}\pp$ associated to the $K$-principal bundle $G\fd \X=G/K$ via the adjoint action of $K$ on $\pp$.

Since $K$ respects the decomposition $\E_\R=\V_\R\oplus\W_\R$, the vector bundle $E_\R$ on $\X$ associated to $G\fd \X$ via the action of $K$ on $\E_\R$ naturally splits as the sum $V_\R\oplus W_\R$. 

Moreover, there exist two natural complex structures belonging to $\SO (S_{|{\W_R}})=\SO(2)$ on the 2-dimensional real vector space $\W_\R$, and only one that we call $I$, such that the orientation of the basis $(Iw,w)$ of $\W_\R$ coincides with the fixed one (for any non-zero $w\in\W_\R$). In the above basis of $\W_\R$, the matrix of $I$ is given by
$$
I=\left(\begin{array}{cc} 0 & 1\\ -1 & 0
\end{array}\right).
$$
The complex structure $I$ then defines a complex structure on the vector bundle $W_\R$ that we still denote by $I$.

Finally, using the identification $T\X\simeq\Hom(W_\R,V_\R)$, we get the complex structure $J$ on $T\X$ : if $X$ is a section of $\Hom(W_\R,V_\R)$, $JX=X\circ I$.
The $G$-invariant Kähler metric $g_\X$ on $\X$ is defined at $o$ by
$$
g_\X(X,Y)=\tr\left(YX\right)=2\,\tr\left(\trans BA\right),\ {\rm if}\ X=\left(\begin{array}{cc} 0 & A\\ \trans A & 0
\end{array}\right),Y=\left(\begin{array}{cc} 0 & B\\ \trans B & 0
\end{array}\right)\in \got p.
$$
The corresponding Kähler form will still be denoted by $\o_\X=g_\X(J.,.)$.  

Next, we consider the complexifications $\V$, $\W$ and $\E=\V\oplus\W$ of $\V_\R$, $\W_\R$ and $\E_\R$ respectively. We extend the  quadratic form $S$ to $\E$ and the complex structure $I$ to $\W$ by $\C$-linearity and still denote them by the same letters.
Let $\W^{1,0}$ (resp. $\W^{0,1}$) be the eigenspace of $I$ corresponding to the eigenvalue $i$ (resp. $-i$). These two eigenspaces also are the two isotropic lines in $\W$ for the quadratic form $S_{|\W}$. Moreover, $\W^{0,1}$ may be identified with $\left({\W^{1,0}}\right)^\star$ by the mean of $S_{|\W}$, and $\V^\star$ may be identified with $\V$ by the mean of $S_{|\V}$.

Let us define $e'_{p+1}=\frac{1}{\sqrt2}(e_{p+1}+ie_{p+2})$ and $e'_{p+2}=\frac{1}{\sqrt2}(e_{p+1}-ie_{p+2})$. In the sequel, we shall use the basis $(e_1,\dots,e_p,e'_{p+1},e'_{p+2})$ of $\E$. The quadratic form $S$ then writes
$$
S=\left(\begin{array}{ccc} I_{p}  & 0 & 0\\ 0 & 0 & -1\\ 0& -1 & 0
\end{array}\right)
$$
and 
$$
\pp=\left\{\left(\begin{array}{ccc} 0 & C & \bar C\\ \trans\bar C &  0 & 0\\ \trans C &0&0 \end{array}\right)\ ,\ C\in M_{p,1}(\C)\right\}.
$$

The complexifications $G_\C$ of $G$ and $K_\C$ of $K$ are respectively  $\SO(p+2,\C)$ and $\SO(p,\C)\times\SO(2,\C)$. 
The Lie algebra $\gg_\C$ of $G_\C$ splits as $\kk_\C\oplus \pp_\C$ where $\kk_\C$ is the Lie algebra of $K_\C$ and 
$$
\pp_\C=\pp\otimes\C=\left\{\left(\begin{array}{ccc} 0 & C & D\\ \trans D &  0 & 0\\ \trans C &0&0 \end{array}\right)\ ,\ C,D\in M_{p,1}(\C)\right\}\,\simeq_\C\, \Hom(\W_\R,\V_\R)\otimes\C=\Hom_\C(\W,\V)
$$

The two eigenspaces of  the extended complex structure $J\otimes{\rm Id}$ acting on $\pp_\C$ are
$$
\pp^{1,0}=\left\{\left(\begin{array}{ccc} 0 & C & 0\\ 0 &  0 & 0\\ \trans C&0&0 \end{array}\right)\ ,\ C\in M_{p,1}(\C)\right\}\,\simeq_\C \Hom_\C(\W^{1,0},\V)\simeq_\C \Hom_\C(\V,\W^{0,1})
$$
and 
$$
\pp^{0,1}=\left\{\left(\begin{array}{ccc} 0 & 0 & D\\ \trans D &0&0\\ 0&  0 & 0 \end{array}\right)\ ,\ D\in M_{p,1}(\C)\right\}\,\simeq_\C \Hom_\C(\W^{0,1},\V)\simeq_\C\Hom_\C(\V,\W^{1,0})~.
$$

The complexified tangent bundle $T^\C\X$ of $\X$ is isomorphic to $G\times_{\Ad K}\pp_\C\simeq(G\times_K K_\C)\times_{\Ad K_\C}\pp_\C$, whereas the holomorphic tangent bundle $T^{1,0}\X$ is isomorphic to $G\times_{\Ad K}\pp^{1,0}\simeq(G\times_K K_\C)\times_{\Ad K_\C}\pp^{1,0}$. There is a natural Hermitian metric on the holomorphic tangent bundle of $\X$ given on $T_o^{1,0}\X\simeq M_{p,1}(\C)$ by $h(C,D)=4\tr\left(\trans\bar DC\right)$.

The holomorphic sectional curvature for the complex line $\la X\ra$ generated by a nonzero $X=\left(\begin{array}{ccc} 0 & C & \bar C\\ \trans\bar C &  0 & 0\\ \trans C &0&0 \end{array}\right)\in T_o\X$ is given by
$K(\la X\ra)=\displaystyle -1+\frac{1}{2}\,\frac{\bigl|\trans CC\bigr|^2}{\left(\trans\bar CC\right)^2}$.
It is clear that $K(\la X\ra)$ is pinched between $-1$ and $-1/2$. The metric $g_\X$ is Einstein and its Ricci curvature tensor is $-\frac{p}{2}\,g_\X$.

\subsection{\texorpdfstring{Toledo invariant and $\SO_0(p,2)$-Higgs bundles}{Toledo invariant and SO(p,2)-Higgs bundles}}\hfill

Let $\rho$ be a representation of a (torsion free) uniform lattice $\G$ of $\SU(m,1)$ into $G=\SO_0(p,2)$, $p\geq 3$. We will assume that $\rho$ is reductive. If it is not, just replace $\rho$ by its semi-simplification $\rho_{ss}$ (see Paragraph~\ref{nonreductive}) in the following to get the result.  

Let $(P_{K_\C},\t)$ be the $G$-principal Higgs bundle on $M$ associated to $\rho$ and $f$ as in section~\ref{higgs} and  
let $E$ be the holomorphic vector bundle on $M$ associated to $P_{K_\C}$ via the action of $K_\C$ on $\E$. The bundle $E$ splits holomorphically as the sum of the rank $p$ subbundle $V=P_{K_\C}\times_{K_\C}\V$ with the rank $2$ subbundle $W=P_{K_\C}\times_{K_\C}\W$. But, in the present situation, we get more structure on $W$, because $K_\C$ also respects the decomposition  $\W=\W^{1,0}\oplus\W^{0,1}$. This implies that if we call $L$ the line bundle $P_{K_\C}\times_{K_\C}\W^{1,0}$, $W$ holomorphically splits as the sum $L\oplus L^{-1}$. Moreover, since $K_\C$ preserves $S_{|\V}$, we can identify $V^\star$ with $V$. In particular, $\deg V=0$.

The Higgs field $\t$ is a holomorphic (1,0)-form taking values in the bundle $P_{K_\C}\times_{\Ad K_\C}\pp_\C\simeq\Hom(L,V)\oplus\Hom(L^{-1},V)$ so that we can write (see also~\cite{BGPG2})
$$
\t=\left(\begin{array}{ccc} 0 & \b & \g\\ \trans\g & 0 & 0\\ \trans\b &0&0\end{array}\right),\mbox{ where }\left\{\begin{array}{l}\b:L\fd V\otimes\Om\\ \g:L^{-1}\fd V\otimes\Om \end{array}\right.
$$ 
The Higgs vector bundle $(E,\t)$ is polystable.  

The bundle $f^\star T^{1,0}\X$ is isomorphic to the bundle $\Hom(L,V)$ and therefore its degree is given by $-p\,\deg L$. We obtain that $\tau(\rho)=\frac{4\pi}{m!}\,\deg L$. Hence the main theorem in this case follows from: 

\begin{theo}\label{ineqso}
$|\deg L| \leq \frac{1}{m}\deg \Om$, that is $|\tau(\rho)|\leq\frac{m+1}{m}\,{\rm Vol}(M)$. In particular, when $m>1$, a representation $\rho:\G\fd\SO_0(p,2)$ is never maximal.
\end{theo}

\begin{demo} 
We shall denote by $\t^2$ the morphism of vector bundles
$$\begin{array}{rcl}
\t^2:T^1\times T^1& \fd & \End (E)\\
(X,Y)&\longmapsto&\t(X)\circ\t(Y)
\end{array},
$$
by $\trans\g\b$ the morphism
$$\begin{array}{rcl}
\trans\g\b:T^1\times T^1& \fd & \End (L)\\
(X,Y)&\longmapsto&\trans\g(X)\circ\b(Y)
\end{array},
$$
etc.

We remark that the system of Hodge bundles obtained after deformation of the Higgs bundle $(E,\t)$ via the $\C^\star$-action on the moduli space is very simple. Indeed, the limiting $(E,\t)$ must verify $\t^n=0$ for some $n$. In particular, as
$$
\t^2=\left(\begin{array}{ccc} \b\trans\g+\g\trans\b & 0 &0\\ 0& \trans\g\b&\trans\g\g\\0&\trans\b\b&\trans\b\g \end{array}\right),
$$
the $(2,2)$ matrix
$$
\left(\begin{array}{cc} \trans\g\b(X,Y)&\trans\g\g(X,Y)\\ \trans\b\b(X,Y)&\trans\b\g(X,Y) \end{array}\right)
=
\left(\begin{array}{cc} \trans\g\b(X,Y)&\trans\g\g(X,Y)\\ \trans\b\b(X,Y)&\trans\b\g(Y,X)\end{array}\right)
$$
(we use $\t^2(X,Y)=\t^2(Y,X)$ and in particular $\trans\b\g(X,Y)={\trans\b}\g(Y,X)$) must be trace free for any $X,Y\in T^1$, which implies $\trans\g\b=0$ (and $\trans\b\g=0$).
Thus
$$
\t^2=\left(\begin{array}{ccc} \b\trans\g+\g\trans\b & 0 &0\\ 0& 0&\trans\g\g\\0&\trans\b\b&0 \end{array}\right).
$$
Suppose now that there exist $X,Y\in T^1$ such that $\trans\g\g(X,Y)\not=0$. For any $X',Y'\in T^1$, $\trans\b\b(X',Y')\trans\g\g(X,Y)=0$ because $\t^n=0$, so we conclude that for every $x\in M$, either $\trans\b\b=0$ or $\trans\g\g=0$ on $T_x^1\times T_x^1$ and then, by holomorphicity, either $\trans\b\b=0$ or $\trans\g\g=0$ on $T^1\times T^1$.


We work on the system of Hodge bundles that we just described. Suppose for example that $\trans\b\b=0$. Then, the sequence
$$
L^{-1}\stackrel{\g}{\fd}\Im\g\stackrel{\trans\g_{|\Im\g}}\fd L\stackrel{\b}{\fd}\Im\b\fd 0
$$ 
defines a Higgs subsheaf of $(E,\t)$.

The bundle $L\otimes T^1$ is semistable and $\Im\b$ is also a Higgs subsheaf of $E$, so we have $\mu(L\otimes T^1)\leq \mu(\Im\b)\leq0$ and then $\deg L\leq\frac{\deg\Om}{m}$ (if $\b=0$, then $\deg L\leq 0$).

Let us consider the maps $\g:L^{-1}\otimes T^1\fd\Im\g$ and $\trans\g_{|\Im\g}:\Im\g\fd L\otimes\Om$. We call $r$ and $r'$ their respective rank. By stability, we have
$$\deg\Im\g\geq r\left(-\deg L+\frac{1}{m}\deg T^1\right)$$
and
$$r'\left(\deg L+\frac{1}{m}{\deg\Om}\right)\geq\deg\Im\trans\g_{|\Im\g}=\deg\Im\g-\deg\Ker\trans\g_{|\Im\g}.
$$
Using the fact that $\deg\Ker\trans\g_{|\Im\g}\leq 0$, we immediately get $\deg L\geq-\frac{\deg\Om}{m}$.
\end{demo}


\end{document}